\documentclass[reqno,12pt]{amsart}

\usepackage{color} 
\usepackage{ifpdf}
\ifpdf
    \usepackage[pdftex]{graphicx}
    \usepackage[pdftex]{hyperref}
    \hypersetup{
        unicode=false,          
        pdftoolbar=true,        
        pdfmenubar=true,        
        pdffitwindow=false,     
        pdfstartview={FitH},    
        pdftitle={MCP Article},      
        pdfauthor={Sara Pollock},   
        pdfsubject={Mathematics},    
        pdfcreator={Sara Pollock},  
        pdfproducer={Sara Pollock}, 
        pdfkeywords={PDE, analysis, numerical analysis}, 
        pdfnewwindow=true,      
        colorlinks=true,        
        linkcolor=red,          
        citecolor=blue,         
        filecolor=magenta,      
        urlcolor=cyan           
    }

\else
    \usepackage{graphicx}
    \usepackage{pstricks}
    
    \newcommand{\href}[2]{#2}
\fi

\usepackage{times}
\usepackage{amsfonts}

\usepackage{amsmath}
\usepackage{amsthm}
\usepackage{amssymb}
\usepackage{amsbsy}
\usepackage{amscd}
\usepackage{extarrows}

\usepackage{enumerate}
\usepackage{verbatim}
\usepackage{subfigure}

\newtheorem{theorem}{Theorem}[section]

\newtheorem{assumption}[theorem]{Assumption}

\newtheorem{example}[theorem]{Example}
\newtheorem{remark}[theorem]{Remark}

\newtheorem{algorithm}[theorem]{Algorithm}
\newtheorem{criteria}[theorem]{Criteria}
\newtheorem{update}[theorem]{Update}

\numberwithin{equation}{section}  

  \newcounter{mnote}
  \let\oldmarginpar\marginpar
    \renewcommand\marginpar[1]{\-\oldmarginpar[\raggedleft\footnotesize #1]%
    {\raggedright\footnotesize #1}}

\newenvironment{enumerateY}
{\begin{list}{{\it(\roman{enumii})}}{
\usecounter{enumii}
\leftmargin 2.5em\topsep 0.5em\itemsep -0.0em\labelwidth 50.0em}}
{\end{list}}

\newenvironment{enumerateX}
{\begin{list}{\arabic{enumi})}{
\usecounter{enumi}
\leftmargin 2.5em\topsep 0.5em\itemsep -0.0em\labelwidth 50.0em}}
{\end{list}}

\definecolor{myblue}{rgb}{0.2,0.2,0.7}
\definecolor{mygreen}{rgb}{0,0.6,0}
\definecolor{mycyan}{rgb}{0,0.6,0.6}
\definecolor{myred}{rgb}{0.9,0.2,0.2}
\definecolor{mymagenta}{rgb}{0.9,0.2,0.9}
\definecolor{mywhite}{rgb}{1.0,1.0,1.0}
\definecolor{myblack}{rgb}{0.0,0.0,0.0}

\renewcommand{\div}{{\operatorname{div}}}
\newcommand{\eps}{\varepsilon}

\newcommand{\PP}{{\mathbb P}}       
\newcommand{\R}{{\mathbb R}}       

\newcommand{\cL}{{\mathcal L}}
\newcommand{\cM}{{\mathcal M}}

\newcommand{\cS}{{\mathcal S}}
\newcommand{\cT}{{\mathcal T}}

\DeclareMathAlphabet{\mathpzc}{OT1}{pzc}{m}{it}

\newcommand{\f}{\frac}

\newcommand{\an}{\text{ and }}

\newcommand{\inn}{\text{ in }}

\newcommand{\tforall}{\text{ for all }}

\newcommand{\bigo}{{\mathcal O}} 

\newcommand{\rest}{\big|}
\newcommand{\Rest}{\Big|}

\newcommand{\grad}{\nabla} 

\newcommand{\goto}{\rightarrow}

\newcommand{\norm}[1]{\ensuremath{\lVert{#1} \rVert}}
\newcommand{\nr}[1]{\norm{#1}} 

\newcommand{\anorm}{\norm{\ \cdot \ }}

\newcommand{\pa}{\partial}

\definecolor{blue}{rgb}{0.2,0.2,0.7}
\definecolor{red}{rgb}{0.7,0.3,0.1}
\definecolor{cyan}{rgb}{0.2,0.5,0.6}
\usepackage{mathtools}
\usepackage{stmaryrd}

\begin{document}

\title[A solver for nonlinear convection diffusion]
      {An improved method for solving quasilinear convection diffusion problems on a coarse mesh}

\author[S. Pollock]{Sara Pollock}
\email{snpolloc@math.tamu.edu}

\address{Department of Mathematics\\
         Texas A\&M University\\ 
         College Station, TX 77843}


\date{\today}



\begin{abstract}
A method is developed for solving quasilinear convection diffusion problems starting on a coarse mesh where the data and solution-dependent coefficients are unresolved, the problem is unstable and approximation properties do not hold.  
The Newton-like iterations of the solver are based on the framework of regularized pseudo-transient continuation where the proposed time integrator is a variation on the Newmark strategy, designed to introduce controllable numerical dissipation and to reduce the fluctuation between the iterates in the coarse mesh regime where the data is rough and the linearized problems are badly conditioned and possibly indefinite.
An algorithm and updated marking strategy is presented to produce a stable sequence of iterates as boundary and internal layers in the data are captured by adaptive mesh partitioning.
The method is suitable for use in an adaptive framework making use of local error indicators to determine mesh refinement and targeted regularization.  
Derivation and  q-linear local convergence of the method is established, and numerical examples demonstrate the theory including the predicted rate of convergence of the iterations.
\end{abstract}

\maketitle




\section{Introduction}
 This paper builds on the framework of~\cite{Pollock14a} and develops a nonlinear solver suitable for use in adaptive methods for quasilinear elliptic problems.  The method is developed to stabilize the linearizations of nonlinear diffusion and convection diffusion problems, especially when there may be steep internal or boundary layers in the problem data.  The sequence of linear problems encountered by a Newton-like method under these circumstances takes the form of  convection diffusion or reaction convection diffusion and the sequence of  approximate solutions is subject to spikes, overshoots and spurious oscillations in the convection-dominated regime.  
 
The present approach  builds on a regularized version of the pseudo-transient continuation method as in for instance~\cite{BaRo80a,KeKe98,CoKeKe02,FoKe05,Deuflhard11,Pollock14a} and the references therein, and  on each mesh refinement seeks a steady-state solution of the nonlinear evolution problem ${\pa}/{\pa t} (Ru) + g(u(t)) = 0$, for positive semidefinite regularization term $R$ in the interest of solving $g(u) = 0$.  In this analysis, further stabilization is  introduced into the time discretization to address the problem of nonphysical oscillations. Time discretizations featuring controllable high-frequency numerical dissipation are well known in the finite element analysis of structures as in for instance~\cite{Newmark59,HiHuTa77,HiHu78a,Hoff88,ChHu93}, and a variation related to these methods referred to here as the $\sigma$-split Newmark update is presently introduced.  This method makes use of the controllable dissipation of the Newmark update and further controls fluctuations between the iterates by freezing a small fraction the linearization about a point with favorable properties.  The $\sigma$-split Newmark method is derived, $q$-linear local convergence with a predictable rate is established, and the method is demonstrated on  three variations of a model problem with steep internal layers.

The goal of the solver in the adaptive setting is to start on a coarse mesh where the problem data is generally not resolved and produce a sequence of transitional states which may not be accurate solutions to the coarse mesh problems, but do allow \emph{a posteriori} error indicators to detect the layers present in the problem data and refine the mesh to be globally fine enough for stability and locally fine enough to achieve accuracy and efficiency.  The sequence of approximate problems can be classified into three phases, the initial phase where the mesh is globally too coarse and quadrature error is high and many features of the problem data remain undetected by the discretization, 
the pre-asymptotic phase where the mesh is fine enough for stability but the problem data is only partially resolved, and the asymptotic regime where the standard existence, uniqueness and approximation properties hold.  
In the coarse mesh regime the solver uses as much stabilization as necessary to produce smooth transitional solutions; in the pre-asymptotic regime the solver adaptively reduces the added stabilization increasing both accuracy and the convergence rate as the data is resolved and the approximate problem becomes less rough; in the asymptotic phase the solver limits to a standard Newton method where the initial guess interpolated from the previous refinement is a good approximation to the solution and the iterations converge quadratically.

The requirements of the adaptive method are that a local error indicator is computed on each refinement and a \textit{reasonable} marking strategy in the sense of~\cite{Siebert10} is employed to determine the next mesh refinement.  A modification of the standard adaptive marking strategy is proposed in which the marked set is determined by two parts: one that refines the elements with the largest indicators, and the other that refines a subset of the coarsest elements of the mesh when the residual from the final Newton-like iteration is significant.  The method reported here improves on~\cite{Pollock14a} both in terms of efficiency and in terms of the strengths of near-singularities it is able to resolve.
  
The remained of the paper is organized as follows.  Section~\ref{sec:continuation} shows the derivation of the pseudo-transient Newmark and $\sigma$-split Newmark methods,~\eqref{alg:Newmark} and~\eqref{alg:sigmaNewmark}. Section~\ref{sec:localconvergence} demonstrates local $q$-linear convergence of the residual for these methods, both with rate $q = 1-1/\gamma$, with $\gamma$ the parameter from the Newmark update.  Section~\ref{sec:algorithm} describes a basic algorithm to implement the solver in an adaptive method, and Section~\ref{sec:numerics} contains the results of numerical experiments using the described adaptive algorithm and~\eqref{alg:sigmaNewmark}.
  
The following notation is used in the remainder of the paper. The function $g(u)$ refers to a specific problem or problem class and $g(x)$ is used in the formal discussion of Newton-like methods. The $n$th iteration subordinate to the $k$th partition $\cT_k$ is denoted $x^n_k$, while $x^n$ is the $n$th iteration on a fixed partition and $x_k$ is the final iteration on the $k$th mesh, taken as the approximate solution on $\cT_k$.  
In defining the weak and bilinear forms in the next section $(u(x) , v(x) ) = \int_\Omega u(x) v(x)\, dx$, and similarly for vector-valued functions. 
 \section{Target problem class}\label{subsec:problem_class}
The nonlinear solver is developed to approximate solutions to the nonlinear problem $g(u)$ = 0,
for polygonal domain $\Omega$ and $g: X \goto Y^\ast$ with $g'(u) \in \cL(X, Y^\ast)$ for real Banach spaces $X \an Y$, particularly where $g(u)$ takes the form of a quasilinear convection diffusion or diffusion problem in divergence form 
\begin{align}\label{eqn:target_strong}
g(u) \coloneqq -\div(\kappa(u) \grad u) + b(u)\cdot \grad u - f(x) = 0 \inn \Omega \subset \R^2, ~u = 0 \text{ on } \pa \Omega.
 \end{align}
or
 \begin{align}\label{eqn:target2_strong}
g(u) \coloneqq -\div(\kappa(u) \grad u) - f(x) = 0 \inn \Omega \subset \R^2, ~u = 0 \text{ on } \pa \Omega.
 \end{align}
Multiplication by a test function and integration by parts over the divergence term yields the  weak form of~\eqref{eqn:target_strong}
\begin{align}\label{eqn:target_weak}
B(u,v) = (\kappa(u) \grad u, \grad v) + (b(u) \cdot \grad u, v)   \tforall V \in Y^\ast,
\end{align}
and linearizing about $u$ at $w \in X$ yields the bilinear form induced by $g'(u)$ 
\begin{align}\label{eqn:target_linear}
B'(u;w,v) & = (\kappa(u) \grad w, \grad v) + (\kappa'(u) w \grad u, \grad v) + (b(u)\cdot \grad w, v)
+ (b'(u) w \grad u, v). 
\end{align}
For $f \in L_2(\Omega) \cap L_\infty(\Omega)$ and $\kappa(u)$ bounded away from zero with 
$\kappa(s), \kappa'(s) \an \kappa''(s)$ bounded as in~\cite{CaRa94}, then~\eqref{eqn:target2_strong} has a unique solution $u \in W^{1,p}(\Omega)$, with $2 < p < \infty$. Both~\eqref{eqn:target_strong} and~\eqref{eqn:target2_strong} fit into the context of~\cite{Xu96} with the assumption that $\kappa(u)$ is bounded and $g'(u): H_0^1(\Omega) \goto H^{-1}(\Omega)$ is an isomorphism, in which case the solution $u$ is an isolated solution. 

The discretized equation is, find $u_h \in X_h$ such that 
$B(u_h,v)=0 \tforall v \in Y_h$ where $X_h \subset X$ and $Y_h \subset Y$ are discrete finite element spaces with respect to triangulation $\cT_h$, where the family of triangulations $\{\cT_h\}_{0 < h < 1}$ is regular and quasi-uniform  in the sense of~\cite{Ciarlet}. Existence, uniqueness and approximation properties of the discrete problems induced by~\eqref{eqn:target_strong} and~\eqref{eqn:target2_strong} are found in~\cite{CaRa94}, ~\cite{Xu96} and the references therein, assuming the mesh is fine enough.
The problem is then to start on a coarse mesh; one that is not fine enough in terms of 
data resolution, stability or approximation properties, and build one that is.
The goal of the solver is to navigate from a coarse to a sufficiently fine mesh where the approximation properties do hold, and to do so by building an both an efficient adaptive mesh and a reasonable initial guess to start the Newton-like iterations on each refinement along the way.

The linearization of~\eqref{eqn:target2_strong} has the form of a  convection-diffusion equation, and ~\eqref{eqn:target_linear} the linearized form of~\eqref{eqn:target_strong} has the form of a  convection reaction diffusion equation
$\alpha(w,v) + \beta(w,v) + \gamma(w,v)$ with
\begin{align}\label{eqn:react_diffus}
\alpha(w,v)&\coloneqq  (\kappa(u) \grad w, \grad v) \\
\beta(w,v) & \coloneqq  (\kappa'(u) w \grad u, \grad v) + (b(u) \cdot \grad w, v) \\
\gamma (w,v) & =  (b'(u) w \grad u, v). 
\end{align}
Using a standard Newton method to solve the nonlinear~\eqref{eqn:target_weak} with~\eqref{eqn:target_linear} generally does not work when there are steep layers present in the problem data, and coarse mesh approximations of the problem are observed to be indefinite and ill-conditioned, consistent with the observations in~\cite{GoSiWa96}.  Many of the problems encountered including formation of spurious spikes, overshoots and instability are symptomatic of the corresponding linear convection-dominated problems~\cite{ElRa01a}.
Techniques form from the finite element analysis of structures~\cite{Newmark59,HiHuTa77,HiHu78a,Hoff88,ChHu93} use numerical integrators featuring high-frequency dissipation to capture lower frequency modes of the solution.  The approach here investigates the use of the Newmark update and a stabilized variation of it as the time-integrators of a pseudo-transient continuation-like method as in~\cite{DePo92,BaRo80a,CoKeKe02,KeKe98, Pollock14a}.

The next section develops the two following methods to improve the convergence of the coarse-mesh iterates in the sequence of linearized problems. 
A positive semidefinite $R$ is used to target specific degrees of freedom for regularization, and the role of $R$ may be seen either from the homotopy perspective as a modification of the path between an initial $u^0$ and  $u^\ast$ that solves $g(u^\ast)=0$, or from the regularization viewpoint as a penalty against certain characteristics of the iterates.  In what follows, $R$ is chosen based on the Laplacian to penalize against high curvature, with degrees of freedom selected for regularization based on an \textit{a posteriori} error indicator.  The selective approach to regularization distinguishes the method presented here and in~\cite{Pollock14a} from the method of pseudo-transient continuation found elsewhere in the literature.
For  $\gamma \ge 1$,  a Newmark update generalizing a backward Euler discretization yields the iteration
\begin{align}\label{alg:Newmark}
\left( \alpha_n R + \gamma\,  g'(x^n) \right) w^n = -g(x^n), \quad x^{n+1} = x^n + w^n.
\end{align}
The proposed $\sigma$-split Newmark update yields the iteration,
\begin{align}\label{alg:sigmaNewmark}
\left( \alpha_n R + \gamma \left\{ (1 - \sigma) g'(\bar x) + \sigma g'(x^n) \right\} \right) w^n = -g(x^n), \quad x^{n+1} = x^n + w^n,
\end{align}
and analysis and demonstration of this last method is the focus of the remainder of the paper.  Local convergence of~\eqref{alg:Newmark} is discussed along the way and local convergence of~\eqref{alg:sigmaNewmark} is established by a perturbation of that result.  The method based on~\eqref{alg:sigmaNewmark} was also tested and found effective on 
a shifted $p$-Laplacian 
\begin{align}\label{eqn:strong_pLaplace}
g(u) = -\div \left\{  \left(  \eps + |\grad u - a|)^{p-2} \right) (\grad u - a)\right\} - f, 
\end{align}
and the related
\begin{align}\label{eqn:strong_GMZ}
g(u) = -\div\left\{ \left( b + \f{1}{  \eps + |\grad u - a|^2} \right) (\grad u - a)\right\} -f, 
\end{align}
problems similar to those investigated in~\cite{BDK12} and~\cite{GaMoZu11}, which reside outside the current target class but still benefit from the regularization techniques described here when starting the adaptive method from a  coarse mesh, especially if the coefficient $\eps <\!<1$.


\section{Continuation methods}\label{sec:continuation}
The homotopy or pseudo-transient continuation method of stabilizing the Newton-like iterations for finding the solution $x^\ast$ of $g(x) = 0$ is developed by discretizing the ODE
\begin{align}\label{eqn:hom001}
\f{\pa (Rx)}{\pa t} + g(x(t)) = 0, \quad x(0) = x^0,
\end{align}
with a positive semidefinite linear functional $R$. In much of the literature, $R$ is taken to be the identity  or a scaled version thereof~\cite{DePo92,KeKe98}, and the references therein; however, the theory is developed for positive-definite functionals other than the identity~\cite{BaRo80a} where it is referred to as the ``s'' method, and positive semi-definite functionals in~\cite{CoKeKe02, Pollock14a}.
This idea can be generalized to discretizing the ODE based on the normal equations formulation of~\eqref{eqn:hom001}
\begin{align}\label{eqn:normal_hom001}
\f{\pa}{\pa t} R^\ast R \, x + g'(x(t))^\ast g(x(t)) = 0, \quad x(0) = x^0,
\end{align}
 with adjoint $R^\ast$ the formal adjoint of $R$ and $g'(x)^\ast$ the formal adjoint of $g'(x)$.
As shown in~\cite{Pollock14a}, the discretization of~\eqref{eqn:normal_hom001} corresponds to the method of Tikhonov regularization.

As discussed in for instance~\cite{BaRo80a,DePo92,KeKe98,CoKeKe02}, letting $x^n$ a finite dimensional approximation to  $x(t_n)$ with $\Delta t_n = t_{n+1}- t_n$, the standard method is to discretize~\eqref{eqn:hom001} by a backward-Euler approximation to ${\pa }/ {\pa t} (R \, x)$ and a linearization of $g'(x^{n+1})$ about $x^n$. In the case of discretizing~\eqref{eqn:normal_hom001}, $g'(x^{n+1})^\ast$ is approximated by $g'(x^n)^\ast$.  
To increase stability of the linearized system, other discretizations of~\eqref{eqn:normal_hom001} are presently considered.  The backward-Euler time discretization is replaced by the more general Newmark update, and the linearization of $g(x^{n+1})$ is split about two distinct points, one the latest iterate and the other yielding a Jacobian with favorable properties.

By the original method of  backward-Euler discretization and first order Taylor expansion about the previous iterate, the resulting Newton-like iteration is given by
\begin{align}\label{eqn:BE_hom002}
\left( \f{1}{\Delta t_n}R  + g'(x^n) \right) w^n = -g(x^n).
\quad x^{n+1} = x^n + w^n. 
\end{align} 
This method increases the stability of the linear system for positive definite and possibly semi-definite $g'(u)$, but runs the risk of shifting the spectrum of the approximate Jacobian towards zero, making the condition significantly worse in the case that $g'(u)$ is indefinite, resulting in large fluctuations between the iterates.  To cope with this situation which occurs in the coarse mesh approximation of quasilinear problems, the iteration based on the normal equations enforcing the shift of the spectrum away from zero is given by
\begin{align}\label{eqn:normal_hom002}
\left( \f{1}{\Delta t_n}R^*R  + g'(x^n)^\ast g'(x^n) \right) w^n = - g'(x^n)g(x^n) , 
\quad x^{n+1} = x^n + w^n. 
\end{align} 
This method is introduced in~\cite{Pollock14a}, where it is  shown that
~\eqref{eqn:normal_hom002} is also found by minimizing the Tikhonov functional $G_\alpha(w)$ for
\begin{align}\label{eqn:normal_hom003}
G_\alpha(w) = \nr{ g'(x^n)w + g(x^n)}_0^2 +  \alpha_n \nr{Rw}_0^2,
\end{align}
with $\alpha_n = 1/\Delta t_n$ and $\anorm_0$ the $L_2$ norm.  As in~\cite{Engl1996}, the necessary and sufficient condition for the minimizer $w$ is $G'_\alpha(w) = 0$, yielding~\eqref{eqn:normal_hom002}.  While successfully increasing the stability of the system, this method remains unsatisfactory due to the increased complexity of solving the system based on the normal equations. A new method is now introduced which  adds stability while preserving the sparsity of the system.  The cost, as shown in Section~\ref{sec:localconvergence} is trading the asymptotically quadratic convergence of~\eqref{eqn:BE_hom002} for $q$-linear convergence.  The proposed algorithm of Section~\ref{sec:algorithm} updates the parameter $\gamma$ of the Newmark method on each mesh refinement as the sequence of linear systems stabilize until the method reduces to the original~\eqref{eqn:BE_hom002} for which quadratic convergence of the error is observed.

\subsection{Time discretization by the Newmark method}\label{subsec:newmark}
The Newmark method~\cite{Newmark59} discretizes the time derivative $\dot u = \pa u/ \pa t $ by
\begin{align}\label{eqn:Newmark001}
x^{n+1} - x^n =  \Delta t_n \left\{ (1 - \gamma) \dot x^n + \gamma \dot x^{n+1} \right\}.
\end{align}
For $\gamma = 1$, ~\eqref{eqn:Newmark001} reduces to the backward Euler discretization described above.
As discussed in~\cite{Newmark59} this time integrator is second order accurate for $\gamma = 1/2$ and introduces nonphysical damping of the high frequency modes for $\gamma > 1/2$ proportional to $\gamma - 1/2$.  More sophisticated $\alpha$ and $\theta$ collocation methods as in~\cite{HiHuTa77,HiHu78a,ChHu93} incorporating the Newmark update are designed to further control numerical dissipation across targeted frequencies.
In the current context the improved resolution and risk of overshoot of these methods designed with two time derivatives in mind is potentially of interest but their implementation in a pseudo-transient continuation setting is beyond the scope of this article. 
For stabilizing the transitional states of the sequence of coarse-mesh problems,  damping of the high-frequency oscillations is the desirable property and more important than a higher order of accuracy.  

Solving~\eqref{eqn:Newmark001} for $\dot x^{n+1}$
\begin{align} \label{eqn:Newmark001a}
\dot x^{n+1} = \f{x^{n+1} - x^n}{\gamma \Delta t_n} - \f{(1 - \gamma)}{\gamma} \dot x^n.
\end{align}
Applying $R \dot x^{n} = - g(x^{n})$ and  $R \dot x^{n+1} = - g(x^{n+1})$ 
 to~\eqref{eqn:Newmark001a}
\begin{align}\label{eqn:normal_newmark002}
\f{1}{\gamma \Delta t}R (x^{n+1} - x^n)  = \f{(1 - \gamma)}{\gamma} g( x^n)  -  g(x^{n+1})
\end{align} 
Linearizing $g(x^{n+1})$ about $x^n$, obtain the Newton-like iteration
\begin{align}\label{eqn:newmark003}
\left( \f{1}{  \Delta t_n}R  +  \gamma\, g'(x^n) \right) w^n = -  g(x^n),
\quad x^{n+1} = x^n + w^n. 
\end{align} 
For $\gamma > 1$ this method comparatively scales down the influence of the right-hand side data; the danger is still allowing and arguably increasing the likeliness of shifting an indefinite left-hand side operator spectrum towards zero.  By the same approach, the iteration based on the normal equations is given by
\begin{align}\label{eqn:normal_newmark003}
\left( \f{1}{  \Delta t_n}R^*R  + \gamma g'(x^n)^\ast g'(x^n) \right) w^n = -  g'(x^n)g(x^n),
\quad x^{n+1} = x^n + w^n. 
\end{align} 
Due to the overdamping effect of setting the parameter $\gamma > 1/2$, the iterations~\eqref{eqn:newmark003} and \eqref{eqn:normal_newmark003} stabilize the solution in some situations where the iteration~\eqref{eqn:normal_hom002} is infeasible.  Similarly to~\eqref{eqn:normal_hom002}-~\eqref{eqn:normal_hom003}, the solution $w^n$ of \eqref{eqn:normal_newmark003} is seen to be the minimizer of the generalized Tikhonov-type functional
\begin{align}\label{eqn:normal_newmark004}
G_\alpha(w) = \nr{  g'(x^n)w + \f 1 \gamma g(x^n)}_0^2 +  \f {\alpha_n}{\gamma} \nr{Rw}_0^2,
\end{align}
with $\alpha_n = 1/\Delta t_n$. 
\subsection{The Newmark update with $\sigma$-splitting}
The Newark discretization of the previous section successfully introduces high-frequency dissipation increasing the stability of the linearized system but suffers the same drawback as the backward-Euler discretization requiring the formulation based on the normal equations in the case of a possibly indefinite Jacobian to control highly unstable sequences of iterates.  The $\sigma$-split Newmark update runs without the use of the normal equations: effectively freezing  a  small fraction of the Jacobian at a point with favorable properties dramatically reduces the fluctuations between iterates in the coarse mesh regime where boundary and internal layers are only partially resolved.

As seen in the derivation, the method can be thought of as splitting the linearization of $g'(x^{n+1})$ about two points, or more precisely as approximating $\Delta x^n$ by a combination of $ \dot x^{n+1}, \dot x^n$ and $\Delta \dot x^n$ where the first is used in the backward-Euler update, the first two are used in the Newmark update with the second adding control of numerical dissipation;  and the third is introduced here to reduce the fluctuation of the Jacobian outside the domain of convergence of the Newton-like iterations.
The method is derived as follows, and is demonstrated in Section~\ref{sec:numerics}  to work without the use of the normal equations in situations where  other methods  including those involving the normal equations are seen to fail. 

Starting with the Newmark update~\eqref{eqn:Newmark001}, each of the time derivative terms on the right is split into two parts
\begin{align}\label{eqn:sigmasplit001}
x^{n+1} - x^n &= \Delta t_n \{  (1 - \gamma)\sigma \dot x^n + (1 - \gamma)(1 - \sigma)  \dot x^n
+  \gamma \sigma\, \dot x^{n+1} + \gamma(1 - \sigma)  \dot x^{n+1}\}
\nonumber \\
&= \Delta t_n \{  (1 -  \gamma \sigma ) \dot x^n + \gamma(1 - \sigma) ( \dot x^{n+1} - \dot x^n)
+\gamma \sigma \, \dot x^{n+1} \}.
\end{align}
Solving for $\dot x^{n+1}$ and applying the relation $R \dot x^{n+1} = - g(x^{n+1})$
\begin{align}\label{eqn:sigmasplit001a}
\f{1}{\gamma \sigma \Delta t_n}R \Delta x^n - \f {(1- \sigma)}{\sigma}R(\dot x^{n+1} - \dot x^n)
= \f{(1 - \gamma \sigma)}{\gamma \sigma} \dot x^n - g(x^{n+1}).
\end{align}
Applying the relation $R \dot x^{n} = - g(x^{n})$ on the left and linearizing $g(x^{n+1})$ about $x^n$ on the right
\begin{align}\label{eqn:sigmasplit001b}
\left( \f{1}{\gamma \sigma \Delta t_n}R  + g'(x^n) \right) \Delta x^n + \f {(1- \sigma)}{\sigma}( g(x^{n+1}) -g( x^n))
= \f{-1}{\gamma \sigma} g(x^{n}).
\end{align}
Linearizing both $g(x^{n+1})  \an g( x^n)$ about fixed $\bar x$, $g(x^{n+1}) -g( x^n) = g'(\bar x) \Delta x^n$. Applying this and multiplying through by $\gamma \sigma$
\begin{align}\label{eqn:sigmasplit001c}
\left( \f{1}{  \Delta t_n} R + \gamma\{  (1 - \sigma) g'(\bar x) + \sigma g'(x^n) \} \right) \Delta x^n 
= - g(x^{n}).
\end{align}
Equation~\eqref{eqn:sigmasplit001c} is the basis for the $\sigma$-split Newmark update.  It contains four parameters: $\alpha_n = 1/\Delta t_n , ~\gamma, \sigma \an \bar x$; while some guidance is provided, a detailed analysis of these parameters will be addressed in subsequent work by the author.  In the current results $\bar x = 0$ is used, and $\bar x = x^0$ the solution from the previous refinement interpolated onto the current mesh has also been observed to work; $\alpha_n$ is chosen as in~\cite{Pollock14a} for the backward-Euler discretization, $\gamma \ge 1$ may be set adaptively, increased to add stability then decreased to speed convergence; and $\sigma \ge \sigma_0\in (0,1)$ is adaptively sent towards one based on the norm of the latest residual on each Newton-like iteration.  It is observed in numerical experiments that $\sigma$ should be close to one, and as seen in Section~\ref{sec:localconvergence} this is a necessary  for the asymptotic $q$-linear convergence at the rate $q = 1-1/\gamma$, in agreement with the rate found using the Newmark discretization without the $\sigma$-splitting.  

\section{Local convergence}\label{sec:localconvergence} 
Local convergence of the residual is established for both algorithms~\eqref{alg:Newmark} based on the Newmark update, and~\eqref{alg:sigmaNewmark} based on the Nemark update with $\sigma$-splitting of the Jacobian.  The second result follows as a perturbation of the first for $\sigma$ close to unity, relying on a weaker set of assumptions.  Both results require the same Lipschitz condition on the Jacobian.  Denote the open ball about $x$ by $B(x,\eps) = \{ y  \, \rest \nr{x-y} < \eps\}$.

\begin{assumption}\label{assume:lipschitz} There exist $\omega_L,\eps > 0$ so that for all $x,y \in B(x^\ast, \eps)$
\[
\nr{g'(x) - g'(y) } \le \omega_L \nr{x- y} \tforall x,y \in B(x^\ast, \eps).
\]
\end{assumption}

\begin{assumption}\label{assume:LC1} (\textit{c.f.} Assumptions 2.2-2.3 of~\cite{CoKeKe02}, and Assumptions 4.1 and 4.4 of ~\cite{Pollock14a}).  There exists $\beta > 0$ so that for positive semidefinite $R$, $\gamma > 1$, and for all  $0 < \alpha_n < \alpha_M$, then for all $x \in B(x^\ast, \eps)$:
\begin{enumerateX}
\item $ \alpha_n R + \gamma \, g'(x)$ is invertible.
\item $\nr{(\alpha_n R + \gamma \, g'(x))^{-1}} \le M_\gamma$.
\item $\nr{ (\alpha_n R)(\alpha_n R + \gamma \, g'(x))^{-1}} \le \f{1}{1 + \beta \gamma /\alpha_n}.$
\end{enumerateX}
\end{assumption}
First, $q$-linear convergence of the residual is established with rate $q = 1 - 1/\gamma$, for $\gamma > 1$. Linear and asymptotic quadratic convergence for the case $\gamma = 1$ is shown in~\cite{Pollock14a}, following from convergence of the error. In the present discussion, convergence of the error is neither used nor shown:  it is observed in numerical experiments that the residual decreases at the predicted rate over iterations  where $\nr{x^{n+1} - x^n}$ may not yet be decreasing; in contrast, for the case $\gamma= 1$ the same quantity displays asymptotically quadratic convergence to zero together with the residual.  The following proof is a variation of Theorem 2.12  of~\cite{Deuflhard11}, where Assumptions~\ref{assume:lipschitz} and~\ref{assume:LC1} replace the affine contravariant Lipschitz condition used in that version of the Newton-Mysovskikh Theorem for the standard Newton method.
\begin{theorem}\label{theorem:localNMc} Let $\alpha_n \le \nr{g(x^n)} \le \alpha_M$ and let Assumptions~\ref{assume:lipschitz} and~\ref{assume:LC1} hold. Define the open set $\cS$ by
\begin{align}\label{eqn:levelsetNMdef}
\cS = \left\{ x \in B(x^\ast, \eps) \, \Rest \, \nr{g(x)} < \f{2 \beta \gamma}{ 2 + \beta \gamma^2 \omega_L M_\gamma^2} \right\},
\end{align}
and suppose $x^n \in \cS$;  then by iteration~\eqref{alg:Newmark}, $x^{n+1} \in \cS$, and the sequence of residuals converges $q$-linearly to zero with asymptotic rate $q = 1-1/\gamma$.
\begin{proof} Let $\Delta x^n = x^{n+1} - x^n$.  Then iteration~\eqref{alg:Newmark} is given by
\begin{align}\label{iter:NewmarkTh}
\left( \alpha_n R + \gamma\,  g'(x^n) \right) \Delta x^n = -g(x^n), \quad x^{n+1} = x^n + \Delta x^n.
\end{align}

By the integral mean-value theorem
\begin{align}\label{eqn:localNMc001}
g(x^{n} + \lambda \Delta x^n) & = g(x^n) + \int_0^\lambda g'(x^n + t \Delta x^n) \Delta x^n \, dt 
\nonumber \\
& =  g(x^n) + \lambda g'(x^n)\Delta x^n + \int_0^\lambda \left( g'(x^n + t \Delta x^n)  - g'(x^n) \right) \Delta x^n \, dt
\nonumber  \\
& = \left(1 - \f \lambda \gamma \right) g(x^n) + \left( \f  \lambda \gamma \right) \alpha_n R \left( \alpha_n R + \gamma g'(x^n) \right)^{-1} g(x^n)
\nonumber  \\
& \quad + \int_0^\lambda \left( g'(x^n + t \Delta x^n)  - g'(x^n) \right) \Delta x^n \, dt 
\end{align}
Applying Assumption~\ref{assume:LC1} (3) to the second term and the Lipschitz condition~\ref{assume:lipschitz} to the third term of~\eqref{eqn:localNMc001}
\begin{align}\label{eqn:localNMc002}
\nr{g(x^n + \lambda \Delta x^n)} & \le \left( 1 - \f \lambda \gamma \right) \nr{g(x^n)} +  \f{\lambda}{\gamma(1 + \beta \gamma /\alpha_n)} \nr{g(x^n)}
+ \f {\lambda^2 \omega_L}{2} \nr{\Delta x^n}^2.
\end{align}
By the iteration~\eqref{iter:NewmarkTh} and Assumption~\ref{assume:LC1} (2), 
$ \nr{\Delta x}^2 \le M_\gamma^2 \nr{g(x^n)}^2$, yielding
\begin{align}\label{eqn:localNMc003}
\nr{g(x^{n} + \lambda \Delta x^n)} & \le \nr{g(x^n)} \left(  \left( 1 - \f \lambda \gamma \right) +  \f{\lambda \alpha_n}{\gamma(\alpha_n + \beta\gamma )} 
+ \f { \omega_L  \lambda^2  M_\gamma^2}{2} \nr{g(x^n)} \right). 
\nonumber \\
\end{align}
By the assumption $\alpha_n \le  \nr{g(x^n)}$, and for $\lambda \in [0,1]$
\begin{align}\label{eqn:localNMc004}
\nr{g(x^{n} + \lambda \Delta x^n)} 
 & \le \nr{g(x^n)} \left(  \left( 1 - \f \lambda \gamma \right) 
+   \left(\f{\lambda}{\gamma} \right) \left( \f{1}{ \beta \gamma} 
+ \f { \gamma \omega_L  M_\gamma^2}{2}\right) \nr{g(x^n)} \right),
\end{align}
for each $\lambda \in [0,1]$ such that $x^n + t \Delta x^n \in \cS$ for all $t \in [0, \lambda]$.
By the logic of ~\cite{Deuflhard11} Theorem 2.12, proceed by contradiction and assume $x^{n+1} \notin \cS$; then there is a smallest $\bar \lambda \in (0,1]$ with $g(x^n + \bar \lambda \Delta x^n) \in \pa \cS$. For that $\bar \lambda$
\begin{align}\label{eqn:localNMc005}
\nr{g(x^n + \bar \lambda \Delta x^n)} & \le  \nr{g(x^n)} \left(  \left( 1 - \f \lambda \gamma \right) 
+   \left(\f{\lambda}{\gamma} \right) \left( \f{1}{ \beta \gamma} 
+ \f { \gamma \omega_L  M_\gamma^2}{2}\right) \nr{g(x^n)} \right) 
\nonumber \\
&< \nr{g(x^n)},
\end{align}
a contradiction.  This shows $x^{n+1} \in \cS$.  To establish the rate of convergence, set $\lambda = 1$ in~\ref{eqn:localNMc004}.
\begin{align}
\nr{g(x^{n+1})}  & \le  \nr{g(x^n)} \left(  \left( 1 - \f 1 \gamma \right) 
+   \left(\f{1}{\gamma} \right) \left( \f{1}{ \beta \gamma} 
+ \f { \gamma \omega_L  M_\gamma^2}{2}\right) \nr{g(x^n)} \right), 
\end{align}
which shows both that $\nr{g(x^{n+1})} < \nr{g(x^n)}$ and the asymptotic $q$-linear rate of $q = 1 - 1 /\gamma$ as $\nr{g(x^n)} \goto 0$ .
\end{proof}
\end{theorem}
\begin{remark}It is observed that the method converges with the predicted rate when the iterates $x^n \notin \cS$ as given by~\eqref{eqn:levelsetNMdef}.  The lapse in the theory appears to be the bound on $\nr{x^n}$ which apparently converges within a smaller set as compared to the residual.
\end{remark}
The asymptotic convergence of the $\sigma$-split algorithm is addressed in the next theorem as a perturbation of Theorem~\ref{theorem:localNMc}.  While it is not necessary in the proof for $g'(\bar x)$ to be positive definite, $\bar x$ is ideally chosen so that $g'(\bar x)$  improves the condition of the approximate Jacobian, allowing Assumptions~\ref{assume:LCsigma} to hold with better constants than~\ref{assume:LC1}. 
\begin{assumption}\label{assume:LCsigma} (\textit{c.f.} Assumption~\ref{assume:LC1}).  There exist $\beta > 0$ so that for $0 < \sigma_0 \le \sigma \le 1$, fixed $\bar x$, positive semidefinite $R$, $\gamma > 1$, and for all  $0 < \alpha_n < \alpha_M$, then for all $x \in B(x^\ast, \eps)$:
\begin{enumerateX}
\item $ \alpha_n R + \gamma\{ (1 - \sigma) g'(\bar x) - \sigma  g'(x) \}$ is invertible.
\item $\nr{(\alpha_n R + \gamma \{(1 - \sigma) g'(\bar x) - \sigma  g'(x) \})^{-1}} \le M_{\gamma\sigma}$.
\item $\nr{ (\alpha_n R)(\alpha_n R + \gamma \{ 1 - \sigma) g(\bar x) - \sigma g'(x) \})^{-1}} \le \f{1}{1 + \beta \gamma /\alpha_n}.$
\end{enumerateX}
\end{assumption}
As with Assumption~\ref{assume:LC1} (3), the third clause agrees with the similar stability bound of~\cite{CoKeKe02} Assumption 2.3, with $\alpha_n$ replaced by $\alpha_n/\gamma$.  Local convergence of the $\sigma$-split Newmark algorithm is established with the same asymptotic rate as in the previous result.
\begin{theorem}\label{theorem:localsigmaNM} Let $\alpha_n \le \nr{g(x^n)} \le \alpha_M$ and let Assumptions~\ref{assume:lipschitz} and~\ref{assume:LCsigma} hold.  Define $\sigma$ by
\begin{align}\label{eqn:sigmadefTh}
\sigma = \max \left\{\sigma_0,\, 1 - \f{\nr{g(x^n)}}{K_0} \right\}, 
\end{align}
for a given $K_0$.  Then there exists $\delta_1 >0$ such that for $x^n$ in  the  open set $\cS$ given by
\begin{align}
\cS = \left\{ x \in B(x^\ast, \eps) \, \Rest \, \nr{g(x)} <\delta_1 \right\},
\end{align}
 and $x^{n+1}$ defined by iteration~\eqref{alg:sigmaNewmark}, it hold that $x^{n+1} \in \cS$, and the sequence of residuals converges $q$-linearly to zero with asymptotic rate $q = 1-1/\gamma$.
\begin{proof} Let $\Delta x^n = x^{n+1} - x^n$.  Then iteration~\eqref{alg:sigmaNewmark} is given by
\begin{align}\label{iter:sigmaNewmarkTh}
\left( \alpha_n R + \gamma \{(1 - \sigma) g'(\bar x) + \sigma  g'(x^n)  \} \right) \Delta x^n = -g(x^n), \quad x^{n+1} = x^n + \Delta x^n.
\end{align}
Much of the proof parallels Theorem~\ref{theorem:localNMc}, and is summarized here.  Starting with the integral mean-value theorem
\begin{align}\label{eqn:localsigmaNMc001}
g(x^{n} + \lambda \Delta x^n) & = g(x^n) + \int_0^\lambda g'(x^n + t \Delta x^n) \Delta x^n \, dt 
\nonumber \\
& =  g(x^n) + \lambda g'(x^n)\Delta x^n 
+ \int_0^\lambda \left( g'(x^n + t \Delta x^n)  - g'(x^n) \right) \Delta x^n \, dt.
\end{align}
By iteration~\eqref{iter:sigmaNewmarkTh}
\begin{align}\label{eqn:localsigmaNMc002}
g'(x^n) \Delta x^n & = -\f{1}{\gamma \sigma} g(x^n) - \f 1 {\gamma \sigma} \alpha_n R \Delta x^n
- \f {1 - \sigma}{\sigma} g'(\bar x) \Delta x^n.
\end{align}
Bounding the  second term on the right of~\eqref{eqn:localsigmaNMc002} by Assumptions~\ref{assume:LCsigma} (3) 
\begin{align}\label{eqn:localsigmaNMc003}
\nr{ \alpha_n R \Delta x^n }&=
\nr{ \alpha_nR  \left( \alpha_n R + \gamma \{(1 - \sigma) g'(\bar x) + \sigma  g'(x^n)  \} \right)^{-1} g(x^n)} 
\nonumber \\ 
& \le \f{1}{1 + \beta\gamma/\alpha_n} \nr{g(x^n)},
\end{align}
and by Assumptions~\ref{assume:lipschitz} and~\ref{assume:LCsigma} (2), then applying $\alpha_n \le \nr{g(x^n)}$ and $\lambda \in [0,1]$
\begin{align}\label{eqn:localsigmaNMc004}
\nr{g(x^{n} + \lambda \Delta x^n)} & \le
\nr{g(x^n)}  \left\{ \left(1 -\f{\lambda}{\gamma \sigma} \right)
+ \f{ (1 - \sigma)}{\sigma} \lambda M_{\gamma \sigma } \nr{g'(\bar x)} 
+\left(\f{\lambda}{\gamma \sigma} \right) \f{1}{(1 + \beta\gamma/\alpha_n)} \right.
\nonumber \\
& \quad  \left. +  \f{\omega_L \lambda^2 M_{\gamma \sigma}^2}{2} \nr{g(x^n)} \right\}
\nonumber \\
& \le \nr{g(x^n)}  \left\{ \left(1 -\f{\lambda}{\gamma \sigma} \right)
+ \f{ (1 - \sigma)}{\sigma} \lambda M_{\gamma \sigma } \nr{g'(\bar x)} \right.
\nonumber \\
& \left. \quad + \left(  \f{\lambda}{\beta \gamma^2 \sigma } \right)  \nr{g(x^n)}
+  \f{ \lambda \omega_L  M_{\gamma \sigma}^2}{2}  \nr{g(x^n)} \right\}. 
\end{align}
Applying~\eqref{eqn:sigmadefTh} to the quantity $(1 - \sigma)/\sigma$  under the assumption $1 - {\nr{g(x^n)}}/{K_0} \ge \sigma_0$ and expanding in orders of    $\nr{g(x^n)}/K_0 <1$,
\begin{align}\label{eqn:localsigmaNMc005}
\f{1 - \sigma}{\sigma} & = \f{\nr{g(x^n)}}{K_0} \f{1}{(1 -  \nr{g(x^n)}/K_0)} 
= \f{ \nr{g(x^n)} }{K_0} \left(  1 + \bigo\left( \f{ \nr{g(x^n)}}{K_0}\right) \right)
\nonumber \\ & = \f{ \nr{g(x^n)}}{K_0} +  \bigo \left(  \f{ \nr{g(x^n)}}{K_0}\right)^2.
\end{align}
Similarly
\begin{align}\label{eqn:localsigmaNMc006}
\f 1 \sigma & = 
 1 + \f{\nr{g(x^n)}}{K_0} + \bigo \left(\f{ \nr{g(x^n)}}{K_0}\right)^2 .
\end{align}
Applying~\eqref{eqn:localsigmaNMc005} and~\eqref{eqn:localsigmaNMc006} to~\eqref{eqn:localsigmaNMc004}, then for any fixed $P >3$
\begin{align}\label{eqn:localsigmaNMc007}
\nr{g(x^{n} + \lambda \Delta x^n)}  &\le
 \nr{g(x^n)} \left\{  \left(1 -\f{\lambda}{\gamma} \right)   + \f{\lambda}{P \gamma}\nr{g(x^n)} \right.
 \nonumber \\
 &\left. \hskip-55pt \times
P \left(    \f{ K_0(\beta \gamma^2 \omega_L M_{\gamma_\sigma}^2/2 + 1)
 +  \beta\gamma( M_{\gamma \sigma} \nr{g(\bar x)} -1 ) } {K_0 \beta \gamma} \right) + \bigo(\nr{g(x^n)}^2) 
 \right\}.
\end{align}
Supposing
\begin{align}\label{eqn:localsigmaNMc008}
 \nr{g(x^n)} \le \delta_0 \coloneqq  \f{K_0 \beta\gamma/P}{K_0(\beta \gamma^2 \omega_L M_{\gamma_\sigma}^2/2 + 1)
 +  \beta\gamma(M_{\gamma \sigma} \nr{g(\bar x)} -1 )},
\end{align}
it follows that
\begin{align}\label{eqn:localsigmaNMc009}
\nr{g(x^{n} + \lambda \Delta x^n)}  &\le \nr{g(x^n)} \left( 1 - \f{(P-1)}{P}\f{\lambda}{ \gamma } + \bigo(\nr{g(x^n)}^2\right),
\end{align}
so there exists  $\delta_1 \in (0, \delta_0]$ for which $\nr{g(x^n)} < \delta_1$ implies
\begin{align}\label{eqn:localsigmaNMc009a}
\nr{g(x^{n} + \lambda \Delta x^n)}  &< \nr{g(x^n)} \left( 1 - \f{(P-2)}{P} \f{\lambda}{\gamma} \right).
\end{align}
The result follows by the same logic as in the proof of Theorem~\ref{theorem:localNMc}.
\end{proof}
\end{theorem}
\begin{remark}It is observed in some problems that the residual $\nr{g(x^n)}$ can decrease at a rate slightly better than the one predicted when $\sigma \ne 1$, whereas for $\sigma=1$ the rate is as predicted.
\end{remark}
\section{Algorithm}\label{sec:algorithm}
The solver using the $\sigma$-Newmark iteration~\eqref{alg:sigmaNewmark} may be implemented in an adaptive method according to the following basic algorithm. To exploit both the stability of method~\eqref{alg:sigmaNewmark} and the quadratic convergence of~\eqref{eqn:BE_hom002}, the parameters $\gamma \an R$ are updated on each refinement, and the parameters $\alpha \an \sigma$ are updated on each iteration of the solver. The main steps of the adaptive method are as follows with the exit criteria and parameter updates specified below.  An example of a targeted regularization term $R$  based on \textit{a posteriori} error indicators is given in Section~\ref{sec:numerics}.
\begin{algorithm}[Basic algorithm for~\eqref{alg:sigmaNewmark}]\label{algorithm:basic} Start with initial $x^0$, $\gamma$, $\sigma_0$. On partition $\cT_k, ~ k = 0, 1, 2, \ldots$
\begin{enumerateX}
\item Compute $R, ~g'(\bar x)$.
\item Set $\alpha_0 = \nr{g(x^0)}$.
\item While Exit criteria~\ref{criteria:exit} are not met on iteration $n-1$:
\begin{enumerateY}
  \item Set $\sigma$ according to~\eqref{eqn:sigmadefTh}.
  \item Solve $(\alpha_n R + \gamma\{ \sigma g'(x^n) + (1-\sigma)g'(\bar x) \}) \Delta x^n = g(x^n)$.
  \item Update $x^{n+1} = x^n + \Delta x^n$.
  \item Update $\alpha_n$.
\end{enumerateY}
  \item Update $\gamma$ for partition $\cT_{k+1}$ according to~\eqref{update:gamma}
\end{enumerateX}
\end{algorithm}

\begin{criteria}[Exit criteria]\label{criteria:exit}  Given a user set tolerance $tol$, an accepted rate of convergence given by $q_{acc} = 1-1/(M \gamma)$ for some constant $M$, \textit{e.g.}, $M=2$, and a maximum number of iterations either chosen as a constant or based on the predicted rate of convergence, exit the solver on partition $\cT_k$ after calculating iterate $x^{n+1}$ when one of the conditions holds.
\begin{enumerateX}
\item $ \nr{g(x^{n+1})} \le tol $.
\item
\begin{enumerateY}
\item  $ \nr{g(x^{n+1})} < \nr{g(x^{n})} $, ~AND
\item $ \nr{g(x^{n})} < \min\{ \nr{g(x^0)}, \nr{g(x_{k-1})}  \} $, ~AND
\item $\nr{g(x^{n+1)}}/\nr{g(x^{n})} < q_{acc}$, ~AND
\item $\nr{g(x^{n+1)}}/\nr{g(x^{n})} > \nr{g(x^{n)}}/\nr{g(x^{n-1})}$.
\end{enumerateY}
\item Maximum number of iterations exceeded.
\end{enumerateX}
\end{criteria}

In terms of the three phases of the solution process in~\cite{Pollock14a}, the final asymptotic regime is characterized by the iterations terminating by Criteria~\ref{criteria:exit} (1); iterations in the pre-asympototic regime terminate with Criteria~\ref{criteria:exit} (2); and iterations terminate with  a mix of Criteria~\ref{criteria:exit} (2) and (3) in the initial phase.  

The second exit criterion~\ref{criteria:exit} (2) merits some explanation, as it allows the iteration to end once reduction of the residual has slowed indicating the iterate has a attained a reasonably stable configuration from which a good prediction about where to refine the mesh may be determined.  The first two clauses require the residual is decreasing, and has decreased below the level given by the previous iterate with at least as much decrease from the initial iteration on the current partition and the residual from the previous partition, if it is well defined.  The third criterion requires error reduction at or close to the predicted rate, and the last that the error reduction between iterates is slowing down. These four criteria together assure the sequence of transitional states is not getting further in the sense of the solver's residual from a converged solution, and prevent situations such as spikes propagating indefinitely across a sequence of partitions.  While spikes, overshoots or other undesirable characteristics may propagate through several refinements, such solutions will eventually not reduce the residual at the specified rate.  When the adaptive mesh is fine enough, these characteristics are observed to smooth out, and otherwise the solver eventually fails and the solution is reset and started again on a finer mesh.

One of the main improvements of~\eqref{alg:sigmaNewmark} over the regularized method of~\cite{Pollock14a} is the generation of more stable sequences of transitional iterates through the pre-asympototic phase. The sequence of pre-asymptotic approximate problems  only partially resolve the data and as internal layers are uncovered over several refinements the condition of these problems is generally bad.  The high-frequency dissipation of the Newmark strategy combined with the stabilization of the iterates produced by the $\sigma$-splitting allow sequences of solutions to propagate through this regime and the mesh to be marked via error indicators leading to the accurate solution of the problem on otherwise coarser meshes than possible if starting the solver on a mesh that is uniformly fine enough to resolve the data.
 In order to make use of the stability and dissipation properties of larger values of $\gamma = \gamma_k$ as well as the asymptotically quadratic convergence if the iterations are stable for $\gamma=1$, the following minimal guidelines are presented, based on the termination criteria above.
\begin{update}[Newmark parameter $\gamma$]\label{update:gamma}. Generally, if the predicted error reduction rate is achieved, $\gamma$ should be reduced; and if the iteration fails, more stability is needed and $\gamma$ should be increased.
\begin{enumerate}
\item Exit criterium~\ref{criteria:exit} (1): If $\gamma_k > 1$ set $1 \le \gamma_{k+1} < \gamma_k.$
\item Exit criterium~\ref{criteria:exit} (2): If $\nr{g(x^{n+1)}}/\nr{g(x^{n})}$ is within tolerance of $1-1/\gamma_k$, set 
 $1 \le \gamma_{k+1} < \gamma_k.$
\item  Exit criterium~\ref{criteria:exit} (3): Set $\gamma_{k+1} > \gamma_k$.
\end{enumerate}
\end{update}

In practice, the residual tends to get reduced below tolerance once $\gamma=1$.  In the results of Section~\ref{sec:numerics}, $\gamma$ is decreased by two when the target rate is achieved, decreased by one if a stable rate lower than predicted is achieved, increased by one when the solver fails, and by two if the maximum number of iteration is exceeded while the iterations are converging below the acceptable rate.  An initial $\gamma_0$ should be chosen large enough to see error reduction on the initial mesh, and not significantly larger.

 The regularization parameter $\alpha_n = 1/\Delta t_n$ is updated by the method described in~\cite{Pollock14a},  repeated here for convenience.  In accordance with the convergence Theorems~\ref{theorem:localNMc} and~\ref{theorem:localsigmaNM}, this strategy assures $\alpha_n \le \nr{g(x^n)}$ so long as the residual is decreasing.
\begin{update}[Tikhonov regularization parameter $\alpha$]\label{update:alpha} Set $\beta_0 = 1$. For $n \ge 1$,
\[
\alpha_{n} = \beta_n \nr{g(x^{n})}, \quad ~\text{with}~ \beta_{n}  = \f {\nr{g(x^{n})}}{\nr{g(x^{n-1})}}.
\]
To reduce rapid fluctuation of $\beta_n$, correct to ensure $\beta_{n-1}/2 \le \beta_n \le 1$ in the case that   $\nr{g(x^{n})}< \nr{g(x^{n-1})}$ and $\beta_n \le 2 \beta_{n-1}$ if $\nr{g(x^{n})}> \nr{g(x^{n-1})}$.
\end{update}

\subsection{Marking strategy}\label{subsec:mark} 
An \textit{a posteriori} error indicator $\eta_T,~ T \in \cT_k$ is assumed available to determine adaptive mesh refinement and as one option for determining a targeted regularization term $R$.  In the current results,  standard  local residual-based element indicators as in for instance~\cite{Stevenson07,FiVe03}  are used for both of these purposes, further described in Section~\ref{sec:numerics}.  Other approaches to solver- and problem-specific regularization and mesh refinement are currently under investigation by the author.  

The goal of the marking strategy is to build a mesh that is globally fine enough for stability, and locally as fine as necessary to achieve the desired accuracy. To improve the efficiency of the method by increasing the stability of the transitional states in the pre-asymptotic phase, the following marking strategy is presented.  Based on exit criteria~\ref{criteria:exit} there are three possible outcomes of the nonlinear solve on refinement $\cT_k$.
\begin{enumerateX}
\item Exit criterium~\ref{criteria:exit} (1): The iterate  $x_k$  has residual  $\nr{g(x_k)} \le  tol$.     
\item Exit criterium~\ref{criteria:exit} (2): The iterate $x_k$ has residual  $\nr{g(x_k)} >  tol$.
\item Exit criterium~\ref{criteria:exit} (3): The solver failed and $x_k$ is reset to zero.
\end{enumerateX}

In the case of~\ref{criteria:exit} (3), the coarsest set of elements in the mesh are refined.  Failure of the solver reflects a \textit{global} rather than a \textit{local} problem.  Starting on a sufficiently coarse mesh, several resets are expected.

In the case of~\ref{criteria:exit} (2) the error indicators $\eta_T$ are computed, and the mesh is refined both according to the elements with the largest local indicators, and according to the coarsest elements with the largest local indicators.  This strategy allows local refinement to take place in order to capture boundary and internal layers to attain eventual accuracy and efficiency while also building the adaptive mesh fine enough to achieve stability.  In meshes that are too coarse to resolve the data and where the approximate problem may have coefficients based on highly inaccurate approximate solutions, nonphysical overshoots often develop in the iterates; while the error indicators in the vicinity of these spikes may be high, refining primarily in these regions exacerbates the problem. As in case~\ref{criteria:exit} (3), a large residual $\nr{g(u_k)}$ predicts a global issue with the mesh. However, valuable information about the near-singularities in the problem data can be predicted from the non-converged iterates, so some local refinement can build a more efficient mesh.  

In the case of~\ref{criteria:exit} (1) the error indicators are computed, and the mesh is refined with respect to the largest local indicators.  Any \textit{reasonable} marking procedure~\cite{Siebert10} for cases~~\ref{criteria:exit} (1)-(2)  may be applied; in particular for case~\ref{criteria:exit} (2) both the element with the largest local indicator must be marked, as well as the coarsest element with the largest local indicator.   In the current results, the D\"orfler marking strategy is used with parameter $\theta$, which for case~\ref{criteria:exit} (2) is split into $\theta = \theta_C + \theta_F$ and the marked sets $\cM_F \subset \cT_k \an \cM_C \subset \cT_k$ are chosen by  sets of least cardinality with
\begin{align}
\sum_{T \in \cM_F}\eta_T^2 \ge \theta_F \sum_{T \in \cT_k} \eta_T^2, \quad
\sum_{T \in \cM_C}\eta_T^2 \ge \theta_C \sum_{T \in \cT_k} \eta_T^2,
\end{align} 
and the marked set $\cM = \cM_F \cup \cM_C$. In case~\ref{criteria:exit} (1), $\theta_F = \theta$ and $\theta_C = 0$. An heuristic choice of $\Phi(\theta) = \theta_C$ is given in Section~\ref{sec:numerics}.

\section{Numerical Examples}\label{sec:numerics} The nonlinear solver~\eqref{alg:sigmaNewmark} is demonstrated on three problems with different structure in their internal layers. In Example~\ref{ex1:convectiondiffusion} results are reported for a model nonlinear convection-diffusion  problem of the form $g(u)$ = 0  which has  a smooth sinusoidal solution. The first variation on the model problem, Example~\ref{ex2:higherfreq} demonstrates the algorithm on the same differential operator with the problem data chosen to generate a higher frequency solution.  Example~\ref{ex3:concentric} shows the results for a related nonlinear diffusion problem with two concentric internal layers.
The adaptive finite element method is implemented using the finite element library FETK~\cite{Holst2001a} and a direct solver is used on each linear system. Both trial and test spaces are taken as the linear finite element space $V_k$ consisting of Lagrange finite elements  $\PP_1$ over partition $\cT_k$ that satisfy the homogeneous Dirichlet boundary conditions.

The local \emph{a posteriori} residual-based indicator, and corresponding jump-based indicator for element $T \in \cT_k$ with $h_T$ the element diameter are given by
\begin{align}
\label{eqn:jump_indicators}
\zeta_T^2(v) &=  h_T \nr{ J_T(v)  }_{L_2(\pa T)}^2 \\
\label{eqn:indicators}
\eta_T^2(v) &= \eta_{\cT_k}^2(v,T) \coloneqq h_T^2 \nr{ g(v)}_{L_2(T)}^2 +  \zeta_T^2(v), 
\end{align}
$J_T(v) \coloneqq \llbracket [\kappa(v) \grad v  \cdot n \rrbracket_{\pa T}$,
with jump
$\llbracket \phi \rrbracket_{\pa T} \coloneqq {\lim_{t \goto 0} \phi(x + t n) - \phi(x - tn)}$, 
 where  $n$ is the appropriate outward normal defined on $\pa T$.  The error estimator on partition $\cT_k$ is given by the $l_2$ sum of indicators $\eta_{\cT_k}^2 = \sum_{T \in \cT_k}\eta_{T}^2$, and similarly for $\zeta_{\cT_k}$. 
 
On each mesh partition $\cT_k$ the penalty matrix $R = R_{k}$ is a localized version of the Laplacian stiffness matrix denoted $R^{global}$.  Let $D^{loc}$ a diagonal matrix of zeros and ones, $v_j$ a vertex of subordinate to partition $\cT_k$, and $\zeta_T = \zeta_T(u^0)$. Then set 
\[
\tilde \psi_k = \sqrt{\text{median}_{T \in \cT_k}(\zeta_T)}, ~\an
\psi_k = \left\{ \begin{array}{ll}
{\tilde \psi_k}^{1/2} , & \text{ if } {\tilde \psi_k} > 1, \\
 \tilde \psi, & \text{ otherwise},
\end{array}
\right.
\]
and 
\begin{align}
D_{jj }^{loc} = \left\{ \begin{array}{ll}
 0, & \text{ if } \zeta_T \le \psi_k \text{ for each element that contains } v_j 
\text{ as a vertex}, \\
1, & \text{ otherwise}.
\end{array} \right.
\end{align}
Then set $R_k = D^{loc}R^{global}D^{loc}$.
\begin{remark} In~\cite{Pollock14a} a similar strategy is applied using the full indicator $\eta_T$ as opposed to the $\zeta_T$ the jump-term alone.  This modification was made to select degrees of freedom to penalize against curvature that show the greatest disparity in curvature as predicted by the jump term.  The role of the selection function $\psi$ is to regularize based on spikes in the indicators leaving the stable dof unregularized to speed convergence.  As the algorithm approaches the asymptotic phase the selection process $\psi$ becomes unimportant because $\alpha_n \le \nr{g(x^n)} \goto 0$.
\end{remark}
In these results the D\"orfler parameters $\theta_C + \theta_F = \theta$ with $\theta = 0.6$, for refinement $k+1$ is  determined according to
\begin{align}
\theta_C = \Phi(\theta) = \theta \left( \f 1 2 + \f 1 \pi \arctan(\nr{g(u_k)}/100 - \pi/2  ) \right).
\end{align}
This form of $\Phi(\theta)$ is chosen to ensure a significant fraction of the marked set is from the coarse mesh until the residual drops below a given threshold.  In the following examples the parameter for setting $\sigma$ by formula~\eqref{eqn:sigmadefTh} are $\sigma \ge \sigma_0 = 0.9$ and $K_0 = 2000$. The tolerance for the nonlinear solver is $tol = 10^{-7}$.
\begin{example}[Nonlinear convection diffusion]\label{ex1:convectiondiffusion} Consider the nonlinear convection diffusion equation on $ \Omega = [0,1] \times [0,1]$,
\begin{align}\label{eqn:ex1problem}
g(u) \coloneqq-\div( \kappa(u) \grad u) + b(u) \cdot \grad u -  f(x,y) =0~ \in \Omega,
~u = 0 \text{ on } \pa \Omega.
\end{align}
The nonlinear diffusion coefficient is given by
\begin{align}\label{eqn:ex_kappadef}
\kappa(s) = k + \f 1{\left( (\epsilon + (s - a)^2\right)},  \quad ~\text{with}~ a = 0.5, \an k = 1.
\end{align}
The nonlinear convection term is given by
\begin{align}\label{eqn:ex_bdef}
b(s) = ( (s-a), (s-a)^2)^T,
\end{align}
and load $f(x,y)$  chosen so the exact solution $u(x,y) = \sin \pi x \sin \pi y$.
\end{example}
Existence of isolated solutions of~\eqref{eqn:ex1problem}-\eqref{eqn:ex_bdef} is discussed in Section~\ref{subsec:problem_class} following
~\cite{Xu96}, assuming the mesh is fine enough.
This problem without the convection term is investigated in~\cite{Pollock14a}, with $\eps = 10^{-3}$; the method here has been observed to run to convergence from a mesh of 144 elements with $\eps = 8\times10^{-5}$, and was tested on smaller values of $\eps$; however the size of the resulting meshes were unwieldy for test-problems.  Data on the iterations in the pre-asymptotic regime is presented here for $\eps = 2\times10^{-4}$ starting on a mesh of 144 elements, with $\gamma = 10$ and run until convergence; and for $\eps = 6 \times 10^{-4}$ run well into the asymptotic regime.  In the pre-asymptotic case  the iterations are  terminated by exit criteria~\ref{criteria:exit} (2), and in most cases where $\gamma$ is reduced it is reduced by one as the iterates are converging at a rate better than predicted. On refinements $20, 21 \an 28$ the final iterates are converging at or close to the predicted rate $q = 1-1/\gamma$, and $\gamma$ is reduced by two.  Table~\ref{tab:ex1_2e4sum} summarizes the computation showing the final residual $\nr{g(u_k)}$, the final ratio $\nr{g(u^{n+1})}/\nr{g(u^n)}$, the final value of $\sigma$ and the Newmark parameter $\gamma$.

\begin{figure}
\centering
\includegraphics[trim=90pt 10pt 40pt 20pt, clip=true,width=0.35\textwidth]{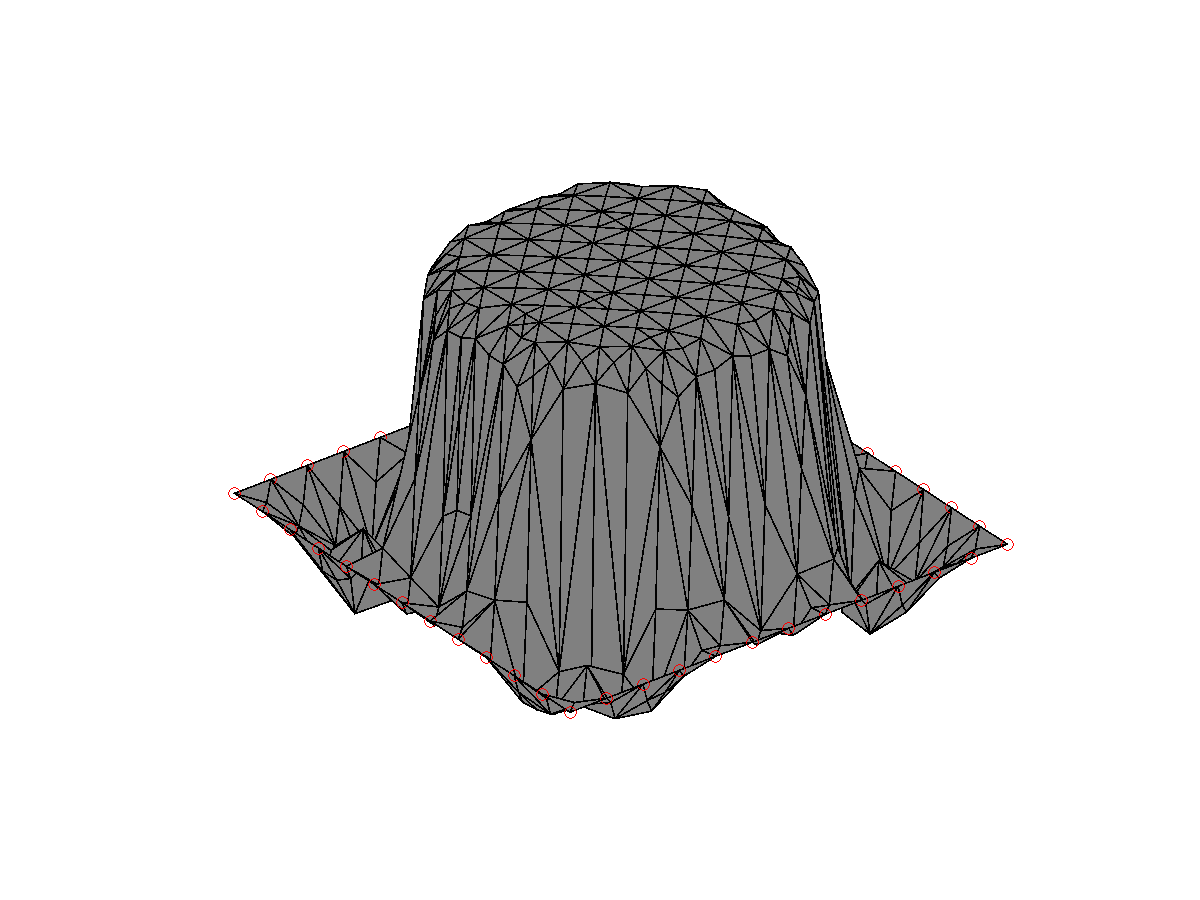}~\hskip-20pt
\includegraphics[trim=90pt 10pt 40pt 20pt, clip=true,width=0.35\textwidth]{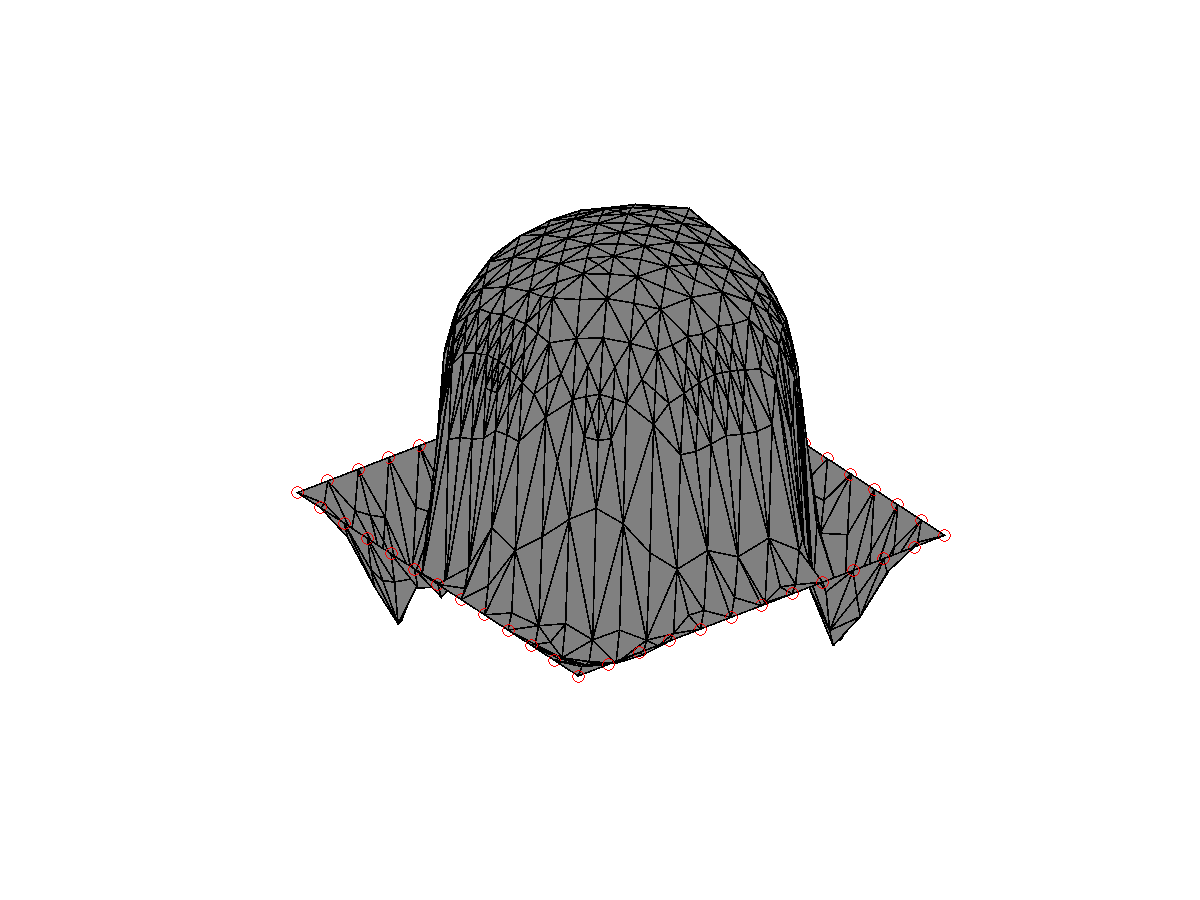}~\hskip-20pt
\includegraphics[trim=90pt 0pt 40pt 20pt, clip=true,width=0.35\textwidth]{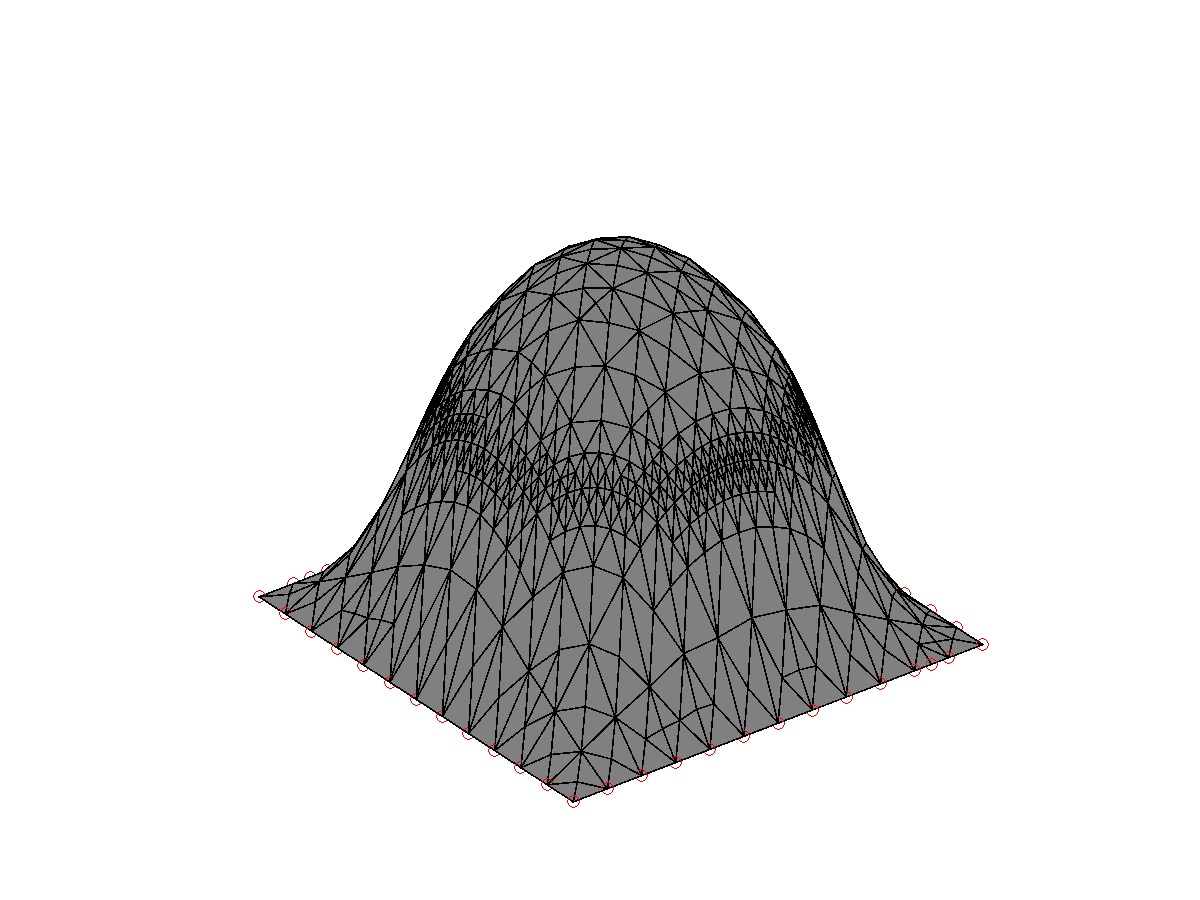}
\caption{Snapshots of the  finite element solution for Example~\ref{ex1:convectiondiffusion} with $\eps = 6\times10^{-4}$ in the pre-asymptotic regime after 6, 12 and 18 adaptive refinements from an initial mesh with 144 elements.}
\label{fig:sol_6e4}
\end{figure}

Of note is the range in magnitude over which the residual shows monotonic and in fact predictable and stable decrease. On level 7, the first partition in the pre-asymptotic regime, the residual is reduced from an initial $\nr{g(u^0)} =  1564.5$, to a final$\nr{g(u^9)} = 871.1$, using $\sigma = \sigma_0 = 0.9$ on each iteration. The final pre-asymptotic solve reduces the residual from $\nr{g(u^0)}= 4.6$, to $\nr{g(u^4)} = 1.4$, with $\sigma >0.996$ and approaching one.   On level 7 for which $\gamma = 20$,  the ratio of consecutive residuals stabilizes well below the predicted rate whereas in level 28, the predicted rate $1-1/\gamma$, with $\gamma = 3$, is achieved.  Most of the pre-asymptotic solves exit after 3-4 iterations, and two of the solves take up to  9; generally the higher number of iterations corresponds to qualitative changes in the iterates such as the smoothing of spikes and changes in curvature illustrated  for instance in the solution snapshots  in Figure~\ref{fig:sol_6e4}.

\begin{table}
	\centering
	\small
\begin{tabular}{c||c|c|c|c}
Level & \raise1pt\hbox{ $\nr{g(u_k)}$}  & \raise4pt\hbox{$\f{\nr{g(u^{n+1})}}{\nr{g(u^n)}}$} & $\sigma_k$ &  $\gamma_k$ \\ 
\hline
\vdots& 			  & 			 &			&		\\
7  	 &871.1	  	  &0.90 	 &0.9			&20		\\
8   	 &786.8 	   	  &0.93 	 &			&20		\\
\noalign{\vskip-18pt}
9   	 &647.0 	  	  &0.94 	 &\raise12pt\hbox{$\vdots$}			&19	 	 \\
10 	 &535.0 	  	  &0.94 	 &			&18		\\
11 	 &445.0 	   	  &0.93 	 &			&17  		\\
12 	 &368.3 	   	  &0.93 	 &			&16 		\\
13 	 &305.2 	  	  &0.93 	 &			&15		\\
\noalign{\vskip-6pt}
14 	 &275.7 	   	  &0.92 	 &$\vdots$      		&14		\\
15 	 &231.3 	  	  &0.92 	 & 0.9		&13	  	\\
16 	 &165.6 	  	  &0.91 	 &0.909		&12 		\\
17 	 &137.7 	  	  &0.90 	 &0.924		&11  		 \\
18 	 &122.1 	  	  &0.90 	 &0.932		&10		\\
19 	 &98.4 	  	  & 0.88 	 &0.944		&9		 \\
20 	 &75.0 	  	  &0.87 	 &0.957		&8		\\
21 	 &59.0 	  	  &0.83 	 &0.964		&6		 \\
22 	 &44.0 	  	  &0.75 	 &0.971		&4		\\
23  	 & 23.5 	  	  &0.67 	 &0.983		&3		\\	
24 	 &14.1 	  	  &0.90 	  &0.990		&3		 \\
25 	 &7.7 	  	  &0.71 	 &0.995		&3		 \\
26 	 &4.8 	   	  &0.68 	 &0.997		&3		 \\
27 	 & 2.5 		  & 0.67 	 &0.998		&3	 	 \\
28 	 &1.4 	  	  & 0.67 	  &0.999		&3		 \\
29 	 &1.9e-08 	 	  &8.1e-03 &1			& 1 		\\
\vdots& 			 & 		&			& \\
\noalign{\vskip3pt}
\end{tabular}
\caption{Summary of data from the pre-asymptotic regime for Example~\ref{ex1:convectiondiffusion} with $\eps = 2\times 10^{-4}$.}
\label{tab:ex1_2e4sum}
\end{table}

\begin{figure}
\centering
\includegraphics[trim=0pt 0pt 0pt 20pt, clip=true,width=0.6\textwidth]{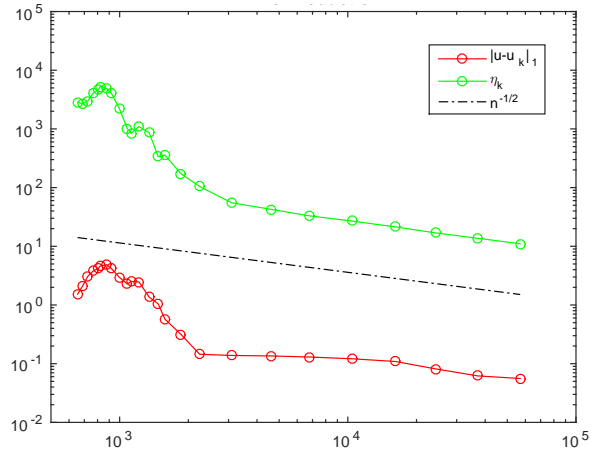}
\caption{$H^1$ error (below) and error estimator (above) against $n^{-1/2}$ where $n$ is the number of elements, for Example~\ref{ex1:convectiondiffusion} with $\eps = 6\times10^{-4}$.}
\label{fig:err_6e4}
\end{figure}

\begin{figure}
\centering
\includegraphics[trim=0pt 0pt 0pt 20pt, clip=true,width=0.3\textwidth]{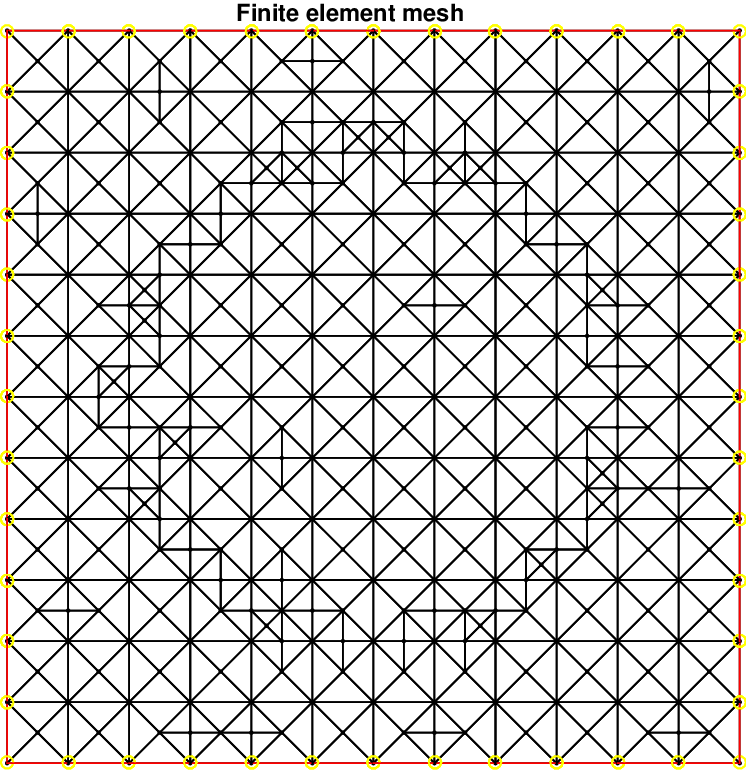}~\hfil~
\includegraphics[trim=0pt 0pt 0pt 20pt, clip=true,width=0.3\textwidth]{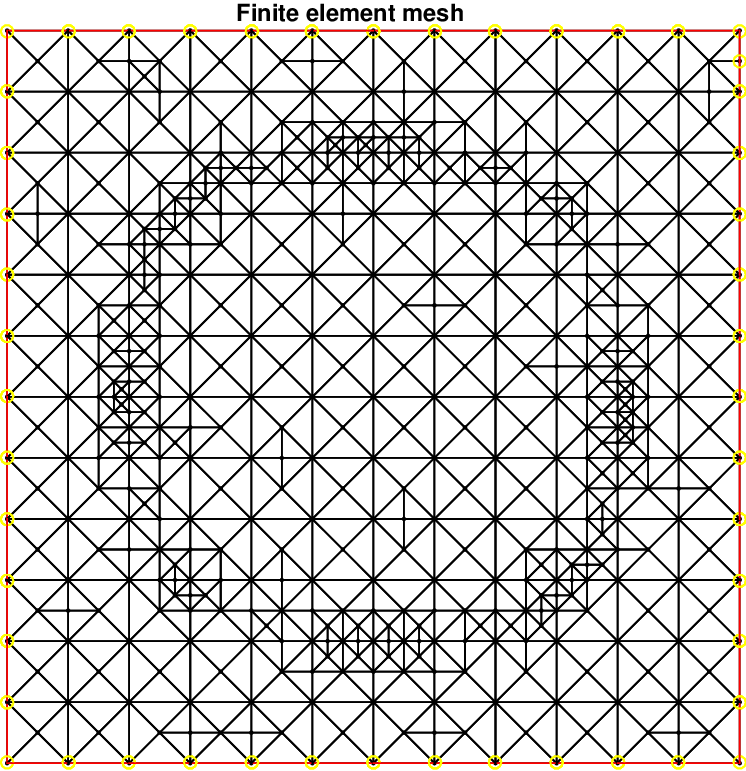}~\hfil~
\includegraphics[trim=0pt 0pt 0pt 20pt, clip=true,width=0.3\textwidth]{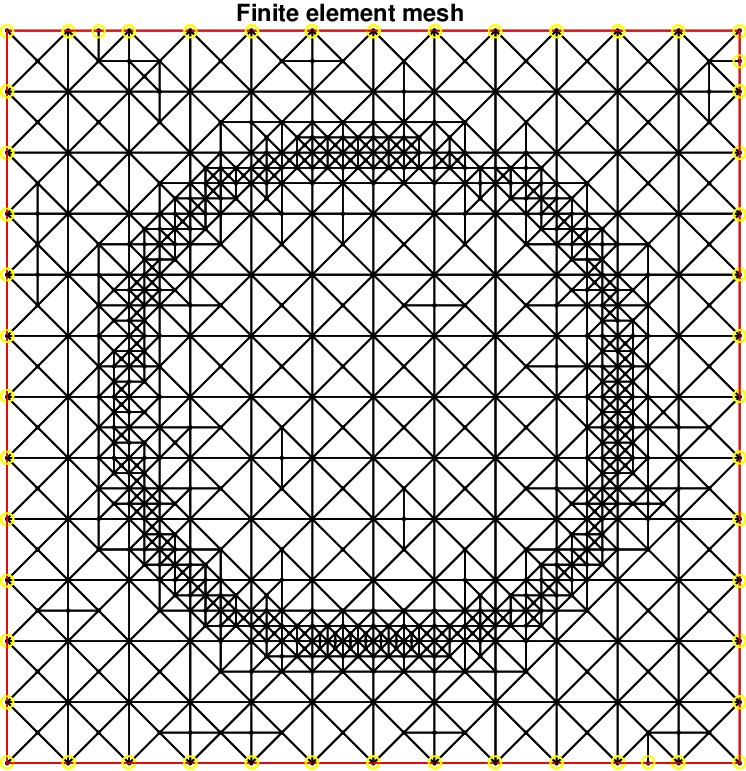}
\caption{Mesh for Example~\ref{ex1:convectiondiffusion} with $\eps = 6\times10^{-4}$ after 6, 12 and 18 adaptive refinements from an initial mesh with 144 elements.}
\label{fig:mesh_6e4}
\end{figure}

\begin{figure}
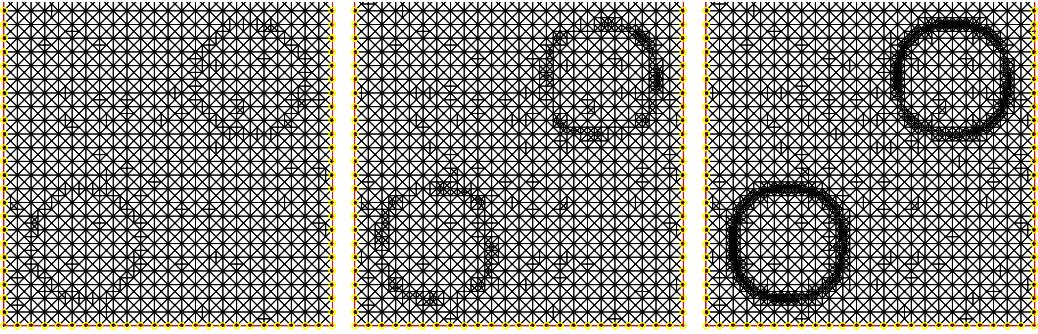

\centering
\includegraphics[trim=0pt 0pt 0pt 20pt, clip=true,width=0.3\textwidth]{mesh_6e4_om2_L8.eps}~\hfil~
\includegraphics[trim=0pt 0pt 0pt 20pt, clip=true,width=0.3\textwidth]{mesh_6e4_om2_L14.eps}~\hfil~
\includegraphics[trim=0pt 0pt 0pt 20pt, clip=true,width=0.3\textwidth]{mesh_6e4_om2_L22.eps}
\caption{Mesh for Example~\ref{ex2:higherfreq} with $\eps = 6\times10^{-4}$ after 8, 14 and 22 adaptive refinements from an initial mesh with 144 elements.}
\label{fig:mesh_6e4_w2}
\end{figure}

\begin{figure}
\centering
\includegraphics[trim=90pt 0pt 40pt 20pt, clip=true,width=0.35\textwidth]{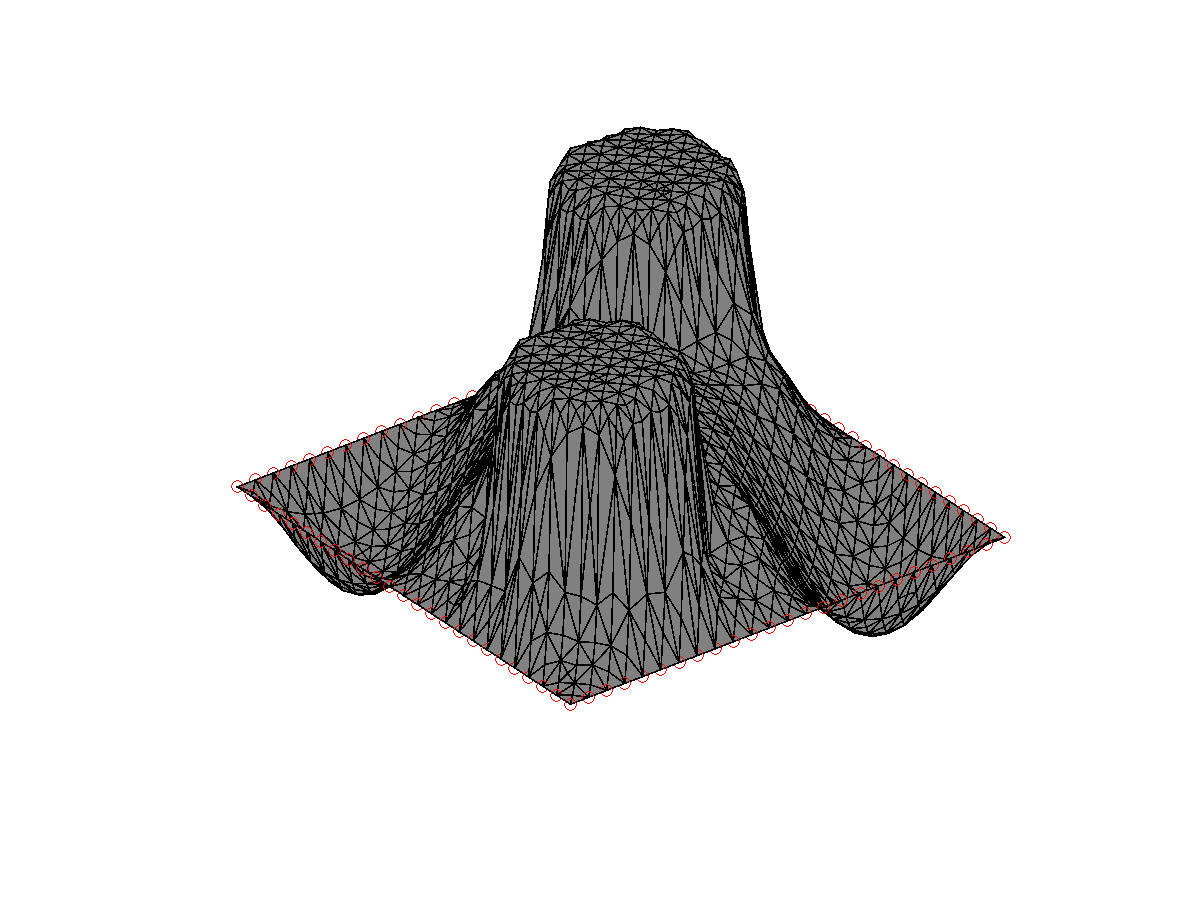}~\hskip-20pt
\includegraphics[trim=90pt 20pt 90pt 20pt, clip=true,width=0.38\textwidth]{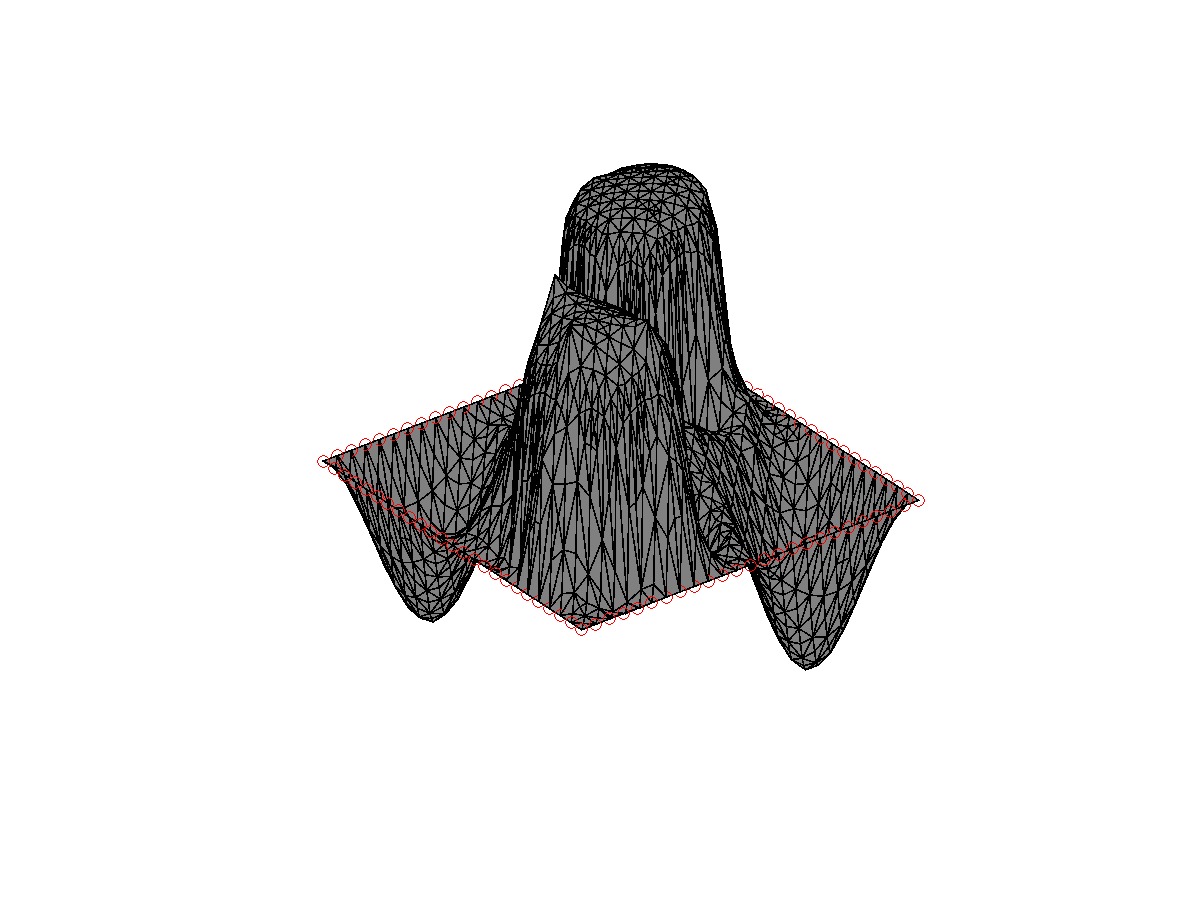}~\hskip-20pt
\includegraphics[trim=140pt 80pt 80pt 20pt, clip=true,width=0.42\textwidth]{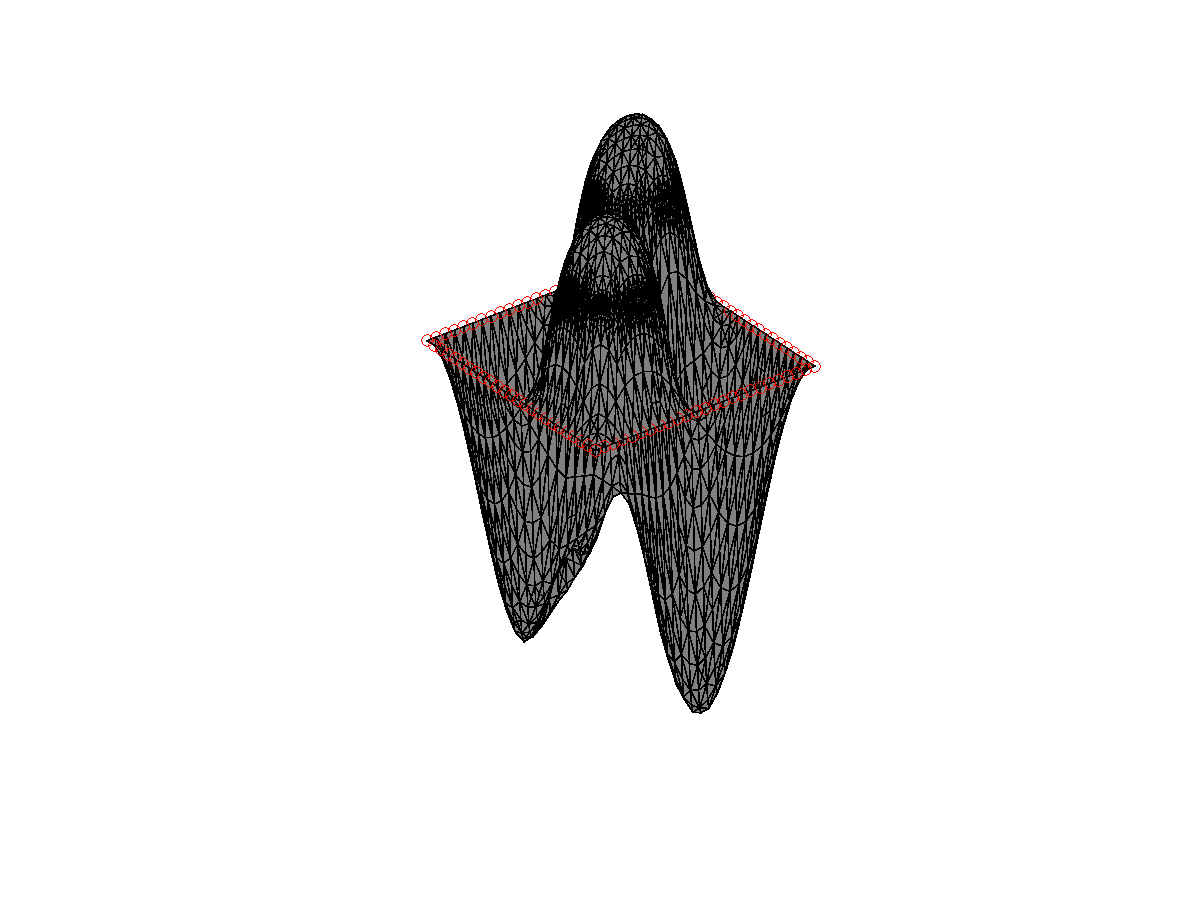}
\caption{Snapshots of the finite element solution for Example~\ref{ex2:higherfreq} with $\eps = 6\times10^{-4}$ in the pre-asymptotic regime after 8, 14 and 22 adaptive refinements from an initial mesh with 144 elements.}
\label{fig:sol_6e4_w2}
\end{figure}

\begin{figure}
\centering
\includegraphics[trim=0pt 0pt 0pt 20pt, clip=true,width=0.6\textwidth]{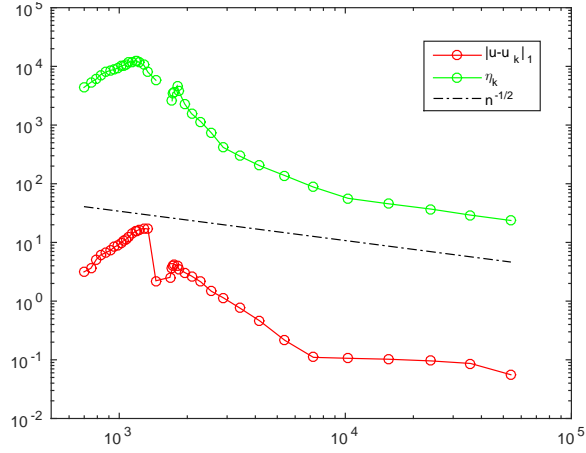}
\caption{$H^1$ error (below) and error estimator (above) against $n^{-1/2}$ where $n$ is the number of elements, for Example~\ref{ex3:concentric} with $\eps = 6\times10^{-4}$.}
\label{fig:err_6e4_a5c8}
\end{figure}

\begin{figure}
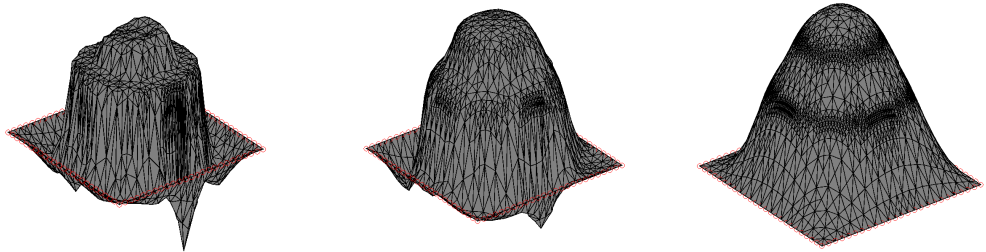

\centering
\includegraphics[trim=90pt 50pt 40pt 20pt, clip=true,width=0.37\textwidth]{sol634_a5c8_L25.eps}~\hskip-20pt
\includegraphics[trim=90pt 20pt 40pt 20pt, clip=true,width=0.35\textwidth]{sol634_a5c8_L30.eps}~\hskip-20pt
\includegraphics[trim=90pt 0pt 40pt 20pt, clip=true,width=0.35\textwidth]{sol634_a5c8_L35.eps}
\caption{Snapshots of the finite element solution for Example~\ref{ex3:concentric} with $\eps = 6\times10^{-4}$ in the pre-asymptotic regime after 25, 30 and 35 adaptive refinements from an initial mesh with 144 elements.}
\label{fig:sol_6e4_a5c8}
\end{figure}

Setting $\eps = 6\times 10^{-4}$, a milder version of the same problem illustrates the reduction in the error from the initial through the pre-asymptotic and into the asymptotic regime which in this case runs from levels $2$ to $19$. The computation is again started on a mesh of 144 elements with $\gamma = 10$.   Figure~\ref{fig:err_6e4} shows a logarithmic plot of the $H^1$ error against the number of elements $n$ compared with $n^{-1/2}$.  The increase in the error and estimator for the first several refinements demonstrates the sequences of partitions capturing the internal layer in the problem data;  the first mesh in Figure~\ref{fig:err_6e4} and corresponding solution in Figure~\ref{fig:mesh_6e4} show the representative behavior as this occurs; the flatness at the top of the solution is characteristic of the penalization by the Laplacian against curvature.  The second mesh of Figure~\ref{fig:err_6e4} and corresponding solution of Figure~\ref{fig:sol_6e4} show the approximate solution as the internal layer is better resolved, but spurious nonphysical oscillations are still apparent in the iterates; this mesh further illustrates the characteristic of the layer getting captured by the mesh, but not uniformly. The final mesh of Figure~\ref{fig:mesh_6e4} and solution of Figure~\ref{fig:sol_6e4} are representative of the end of the asymptotic regime where the mesh now refines to capture the layer with increasingly uniform resolution and the solution has the overall correct shape.  Within several iteration the approximate solution looks like a sinusoid and the Newton-like iterations solve to tolerance.
It is noted in Figure~\ref{fig:err_6e4} that the error is stable for the first several iterations in the asymptotic regime while the error estimator reduces at a steady rate; presumably this occurs due to the presence of another nearby solution of~\eqref{eqn:ex1problem}. A similar phenomenon is observed in~\cite{GaMoZu11} in examples where the monotonicity is violated and the solution is known only to be locally unique. 
\begin{example}[Higher frequency solution]\label{ex2:higherfreq} A variation of Example~\ref{ex1:convectiondiffusion} is considered where the problem is given by~\eqref{eqn:ex1problem}-\eqref{eqn:ex_bdef},  with $\eps = 6 \times 10^{-4}$, and the load $f(x,y)$ chosen so the exact solution is $u = \sin(2 \pi x) \sin(2 \pi y)$.
\end{example}
Local uniqueness of the solution follows from~\cite{Xu96}, assuming the mesh is fine enough.
This example features two disjoint internal layers in the problem data, one about each positive peak of the sinusoid.  The decrease in the residual in the pre-asymptotic regime closely resembles the data presented in Table~\ref{tab:ex1_2e4sum} for Example~\ref{ex1:convectiondiffusion}. For $\eps = 6 \times 10^{-4}$, the pre-asymptotic phase for this example runs from levels $4$ to $26$ as opposed to $2$ to $19$ for the sinusoid with a single peak; and the mesh in Example~\ref{ex2:higherfreq} maintains a maximum meshsize of $h = 8.6\times 10^{-4}$, whereas in Example~\ref{ex1:convectiondiffusion}, the asymptotic meshsize is $h = 3.5 \times 10^{-3}$,  indicating two more coarse mesh refinements are taken to stabilize the iterates.  This is to be expected due to the decrease in width of the corresponding layers.  Mesh partitions and snapshots of the pre-asymptotic iterates are shown in Figures~\ref{fig:mesh_6e4_w2} and~\ref{fig:sol_6e4_w2}.  
The characteristic flatness of the peaks due to the penalization against curvature is again observed in the first two iterates, shown respectively at the $8$th and $12$th adaptive refinements, and resolved by the $22$nd, as seen on the right.  The iterates show a qualitatively different behavior than the nonphysical oscillations of Example~\ref{ex1:convectiondiffusion}; here,  the downward peaks which should have the same magnitude as the upward peaks start as shallow and extend to their full depth as the refinements progress.
\begin{example}[Concentric layers]\label{ex3:concentric} Consider the nonlinear diffusion problem on $ \Omega = [0,1] \times [0,1]$ with two concentric internal layers.
\begin{align}\label{eqn:ex3problem}
g(u) \coloneqq-\div( \kappa(u) \grad u) -  f(x,y) =0~ \in \Omega,
~u = 0 \text{ on } \pa \Omega,
\end{align}
with the nonlinear diffusion coefficient  given by
\begin{align}\label{eqn:ex3_kappadef}
\kappa(s) = k + \f 1{\left( (\epsilon + (s - a)^2\right)} +  \f 1{\left( (\epsilon + (s - c)^2\right)} ,
  \quad ~\text{with}~ a = 0.5, ~c = 0.8 ~\an k = 1.
\end{align}
The load $f(x,y)$ is chosen so the exact solution $u(x,y) = \sin (\pi x) \sin (\pi y)$.
\end{example}
Local uniqueness of the solution again follows from~\cite{Xu96}, assuming the mesh is fine enough.  Similarly to the error decrease in the first example, Figure~\ref{fig:err_6e4_a5c8} shows a temporary leveling off of the error while the  estimator decreases at a steady rate at the beginning of the asymptotic phase indicating the presence of a nearby solution.  The discontinuity in the error estimator and corresponding jump in the error is due to the final coarse mesh refinement where the solution is reset at adaptive level 22, after which the mesh maintains a maximum meshsize of $h = 1.7\times 10^{-3}$, half the meshsize of Example~\ref{ex1:convectiondiffusion}. Once in the pre-asymptotic regime  the internal layers are progressively resolved and the residual drops below tolerance converging at a quadratic rate when $\gamma = 1$, at adaptive level 37.  Due to the aforementioned flattening of the iterates induced by the solver's stabilization, the approximate solutions gradually increase in height as the regularization parameters are reduced; as such, the inner internal layer centered at $u = 0.8$, is uncovered later than the outer layer centered at $u = 0.5$.  This phenomenon is illustrated in the adaptive meshes at levels $25, ~30 \an 35$ shown in Figure~\ref{fig:mesh_6e4_a5c8}.  The corresponding iterates shown in Figure~\ref{fig:sol_6e4_a5c8} display the characteristic spikes between the outer layer and the boundary similar to those in Example~\ref{ex1:convectiondiffusion}, and in this case the peak of the sinusoid progressively resolves the curvature in two flatter regions.  
\begin{figure}
\centering
\includegraphics[trim=0pt 0pt 0pt 20pt, clip=true,width=0.3\textwidth]{mesh634_a5c8_L25.eps}~\hfil~
\includegraphics[trim=0pt 0pt 0pt 20pt, clip=true,width=0.3\textwidth]{mesh634_a5c8_L30.eps}~\hfil~
\includegraphics[trim=0pt 0pt 0pt 20pt, clip=true,width=0.3\textwidth]{mesh634_a5c8_L35.eps}
\caption{Mesh for Example~\ref{ex3:concentric} with $\eps = 6\times10^{-4}$ after 25, 30 and 35 adaptive refinements from an initial mesh with 144 elements.}
\label{fig:mesh_6e4_a5c8}
\end{figure}
\section{Conclusion}\label{sec:conclusion}
The $\sigma$-split Newmark method of regularized pseudo-transient continuation is introduced to solve the nonlinear convection diffusion problem $g(u)= 0$ starting from a coarse mesh where the problem data is not yet resolved.  The method is derived and local $q$-linear convergence in agreement with the observed  rate is established.  An algorithm is presented to fit the solver into a standard adaptive method, and is demonstrated on three variations of a problem with steep internal layers.  The method including the adaptive update of the solver's parameters on each refinement builds on the framework presented in~\cite{Pollock14a},  and is designed to effectively stabilize linearizations over  rough problem data where spikes, overshoots and spurious oscillations otherwise prevent the sequence of transitional solutions from approaching an accurate approximation. The results here are an improvement  first in terms of efficiency  of the new solver which functions without the use of the normal equations; and second by an improved set of exit criteria for the solver which works together with an updated marking strategy to improve the stability of the method through the pre-asymptotic coarse mesh regime.  Ultimately this yields a more efficient route to the asymptotic regime where the problem data is resolved and approximation properties hold. Future work will address targeted local regularization and parameter selection of the nonlinear solver, and  more general types of nonlinearites.

\section{Acknowledgments}
   \label{sec:ack}
The author would like to thank William Rundell and Alan Demlow for providing input on a draft of this work.

\bibliographystyle{siam}
\bibliography{refsTRN}



\end{document}